\documentclass[onefignum,onetabnum]{siamart171218}
\usepackage{bbm}
\usepackage{subfig}
\usepackage{amsfonts}
\usepackage{graphicx}
\usepackage{epstopdf}
\usepackage{algorithmic}
\usepackage{ulem}
\DeclareMathOperator*{\argmin}{arg\,min}

\ifpdf
  \DeclareGraphicsExtensions{.eps,.pdf,.png,.jpg}
\else
  \DeclareGraphicsExtensions{.eps}
\fi

\newsiamremark{remark}{Remark}
\newsiamremark{hypothesis}{Hypothesis}
\crefname{hypothesis}{Hypothesis}{Hypotheses}
\newsiamthm{claim}{Claim}
\crefname{claim}{Claim}{Claim}

\newcommand{\op}{\mathcal{L}}
\headers{GF}{GF}

\title{High order, semi-implicit, energy stable schemes for gradient flows\thanks{\textbf{Funding:} Alexander Zaitzeff and Selim Esedoglu gratefully acknowledge support from the NSF grant DMS-1719727. Krishnna Garikipati acknowledges NSF grant DMREF-1729166.}}

\author{Alexander Zaitzeff\thanks{Department of Mathematics, University of Michigan, Ann Arbor, MI 48109, USA
  (\email{azaitzef@umich.edu}, \email{esedoglu@umich.edu}).}
\and Selim Esedo\=glu\footnotemark[2]
\and Krishna Garikipati\thanks{Departments of Mechanical Engineering, and Mathematics, Michigan Institute for Computational Discovery \& Engineering, University of Michigan, Ann Arbor, MI 48109, USA (\email{krishna@umich.edu}).} }

\usepackage{amsopn}

\ifpdf
\hypersetup{
  pdftitle={Semi-Implicit Schemes for General Gradient Flows},
  pdfauthor={ZKE}
}
\fi

\begin{document}

\maketitle


\begin{abstract}
We introduce a class of high order accurate, semi-implicit Runge-Kutta schemes in the general setting of evolution equations that arise as gradient flow for a cost function, possibly with respect to an inner product that depends on the solution, and we establish their energy stability.
This class includes as a special case high order, unconditionally stable schemes obtained via convexity splitting.
The new schemes are demonstrated on a variety of gradient flows, including partial differential equations that are gradient flow with respect to the Wasserstein (mass transport) distance.
\end{abstract}



\section{Introduction}
We are concerned with numerical schemes for evolution equations that arise as gradient flow (steepest descent) for an energy $f:H \rightarrow \mathbb{R}$, where $H$ is a Hilbert space with inner product $\langle\cdot,\cdot\rangle$:
\begin{equation}
    \label{eq:de}
    u'=-\nabla_H E(u).
\end{equation}
Additionally, we will study gradient flows with a solution dependent inner product:
\begin{equation}
    \label{eq:Lde}
    u'=-\op(u)\nabla_H E(u)
\end{equation}
where $\op(u)$ is a positive definite operator that depends on $u$.

Equations \cref{eq:de} and \cref{eq:Lde} may represent (scalar or vectorial) ordinary or partial differential equations.
One property of \cref{eq:de} and \cref{eq:Lde} is dissipation: $\frac{d}{dt} E(u) \leq 0$, to see this
\[
   \frac{d}{dt} E(u)=\langle \nabla_H E(u),u' \rangle = -\langle \nabla_H E(u),\op(u) \nabla_H E(u) \rangle \leq 0.
\]
In \cite{alexander2019variational} the authors focused on unconditionally stable numerical methods to solve \cref{eq:de}. In this paper,
our focus is on semi-implicit methods that come with a rigorously established energy stability property inherited from the above expression of dissipation. Specifically let
\begin{equation}
\label{eq:additive}
E(u) = E_1(u)+E_2(u)
\end{equation}
where in our numerical implementation we will handle $E_1$ implicitly and $E_2$ explicitly.
We allow any choice for $E_1$ and $E_2$ as long as (\ref{eq:additive}) is satisfied.
Our numerical methods will guarantee that when the time step is less than a constant depending only on $E_2$ the following numerical dissipation property will hold: 
\begin{equation}
\label{eq:energymin}
    E(u_{n+1})\leq E(u_n)
\end{equation}
where $u_n$ denotes the approximation to the solution at the $n$-th time step. 

A basic semi-implicit scheme for the abstract equation \cref{eq:de}, with time step size $k>0$, reads
\begin{equation}
\label{eq:si}
\frac{u_{n+1} - u_n}{k} = -\nabla_H E_1(u_{n+1})-\nabla_H E_2(u_{n}).
\end{equation}
Let $L_2(u,u_n)$ be the linearization of $E_2$ around $u_n$ so 
\begin{equation*}
    L_2(u,u_n)=E_2(u_n)+\langle\nabla_H E_2(u_n),u-u_n\rangle
\end{equation*}
Then \cref{eq:si} is the Euler-Lagrange equation for the optimization problem
\begin{equation}
\label{eq:minmov}
u_{n+1}=\argmin_u E_1(u)+L_2(u,u_n)+\frac{1}{2k}\left\|u-u_n\right\|^2
\end{equation}
where $\| \cdot \|^2 = \langle \cdot , \cdot \rangle$. For \begin{equation}
    \label{eq:lmda}
    \Lambda=\max\{0,\max_{u,\left\|v\right\|=1}D^2E_2(u)\big(v ,v\big)\}
\end{equation}
where by $D^2E_2(u)\big(v ,w\big)$ we mean
$\left. \frac{d^2}{d\epsilon_2d\epsilon_1}E_2(u+\epsilon_1 v+\epsilon_2 w) \right|_{\epsilon_1=\epsilon_2=0}.$
We have
\begin{equation}
\label{eq:bound} 
E_2(u)\leq L_2(u,p)+\frac{\Lambda}{2}\left\|u-p\right\|^2
\end{equation} 
for any $u$ and $p$.
It follows that when $k \leq \frac{1}{\Lambda}$ 
\begin{align*}
E(u_{n+1})&=E_1(u_{n+1})+E_2(u_{n+1})\leq E_1(u_{n+1})+L_2(u_{n+1},u_n)+\frac{1}{2k}\left\|u_{n+1}-u_n\right\|^2\\&\leq E_1(u_{n})+L_2(u_{n},u_n)+\frac{1}{2k}\left\|u_n-u_n\right\|^2=E_1(u_{n})+E_2(u_{n})=E(u_{n})
\end{align*}
so that scheme \cref{eq:si} is stable under this condition on the time step $k$, provided that optimization problem \cref{eq:minmov} can be solved.


Our new class of methods has no assumption (e.g. convexity, concavity) on the components $E_1$ and $E_2$ of the energy that are treated implicitly and explicitly, respectively. 
They have high (at least up to third) order accuracy.
Our stability results are conditional, but revert to unconditional stability when $E_1$ and $E_2$ have appropriate convexity properties, and contain as a special case previous unconditional stability results for high order, convexity splitting type schemes.
In particular, a previous paper on higher order ARK IMEX energy stable schemes \cite{shin2017unconditionally} studies, in the spirit of convexity splitting \cite{eyre1998unconditionally}, formulations that break up the energy into convex and concave parts, and treat the convex part implicitly and the concave part explicitly. We know of no other work on ARK IMEX methods that considers energy stability.

Another novelty of the present paper is extending these high order, stable, implicit and semi-implicit methods for general gradient flows to solve \cref{eq:Lde}, the case when the inner product is solution dependent.
There are certainly existing stable methods for \cref{eq:Lde} on a case by case basis, for example for the Cahn-Hillard equation with degenerate mobility \cite{chen2019positivity,han2015second} and the porous medium equation \cite{del2018robust, duan2019numerical, westdickenberg2010variational}.
To our knowledge, this is the first time energy stable, high order schemes for general gradient flows with solution dependent inner products have been developed.

\noindent 

The rest of the paper is organized as follows:
\begin{itemize}
\item \Cref{sec:gf} presents the conditions for energy stability and constructs schemes that satisfy them.
\item In \cref{sec:ex}, we state the consistency equations for the ARK IMEX schemes for solving gradient flows \cref{eq:de} and give 2nd and 3rd order examples.
\item \Cref{sec:SDIP} gives 2nd and 3rd order methods for solving gradient flows with solution dependent inner product \cref{eq:Lde} and provides consistency calculations.
\item In \cref{sec:nr}, we present numerical convergence studies several of well-known partial differential equations that are gradient flows, including with respect to Wasserstein metrics.
\end{itemize}
The code for \cref{sec:nr} is publicly available, and can be found at \url{https://github.com/AZaitzeff/SIgradflow}.

\section{Stability of Our New Schemes}
\label{sec:gf}
In this section, we formulate a wide class of numerical schemes that are energy stable by construction. The first of these schemes are Implicit-Explicit Additive Runge-Kutta (ARK IMEX) schemes, but  we will write them in variational form in order to prove energy stability more easily.
The variational formulation of our $M$-stage ARK IMEX scheme is:
\begin{enumerate}
    \item Set $U_0 = u_n$.
    \item For $m=1,\ldots,M$:
\begin{equation}
\label{eq:ms}
U_m=\argmin_u \bigg( E_1(u)+\sum^{m-1}_{i=0} \theta_{m,i} L_2(u,U_i)+\sum^{m-1}_{i=0}\frac{\gamma_{m,i}}{2k}\left\|u-U_i\right\|^2 \bigg).
\end{equation}
where
\begin{equation}
\label{eq:L2}
    L_2(u,p)=E_2(p)+\langle\nabla_H E_2(p),u-p\rangle
\end{equation}
\item Set $u_{n+1}=U_M$.
\end{enumerate}
\medskip

\noindent The schemes for approximating a gradient flow with respect to a {\it solution dependent inner product} will be a series of embedded ARK IMEX methods.
The inner product will be fixed for each ARK IMEX step, allowing the stability results of this section to apply.
Now we establish quite broad conditions on the coefficients $\gamma_{m,i}$ $\theta_{m,i}$ that ensure conditional energy dissipation \cref{eq:energymin}. 
Before we state and prove the conditions in generality, consider the following two-stage special case of scheme \cref{eq:ms}:
\small
\begin{align}
\label{eq:twostage1}
U_1=\argmin_u\bigg( E_1(u) &+L_2(u,u_n)+\frac{\gamma_{1,0}}{2k}\left\|u-u_n\right\|^2\bigg)\\
\label{eq:twostage2}
\begin{split}
u_{n+1}=\argmin_u \bigg(E_2(u) &+\theta_{2,1} L_2(u,U_1)+\theta_{2,0} L_2(u,u_n)\\&+\frac{\gamma_{2,0}}{2k}\left\|u-u_n\right\|^2+\frac{\gamma_{2,1}}{2k}\left\|u-U_1\right\|^2\bigg)
\end{split}
\end{align}
\normalsize
Let $\Lambda=\max\{0,\max_{x,\left\|v\right\|=1} D^2E_2(x)\big(v ,v\big)\}$. Note that this implies
\begin{equation}
\label{eq:twostage3} 
E_2(u)\leq L_2(u,p)+\frac{\Lambda}{2}\left\|u-p\right\|^2
\end{equation}
for any $u$ and $p$. Also note that $L_2(u,u)=E_2(u)$.
Impose the conditions
\small
\begin{align}
\label{eq:twostage4} 
\begin{split}
&\gamma_{1,0}-k\Lambda-\frac{(\gamma_{2,0}-k\Lambda\theta_{2,0})^2}{(\gamma_{2,0}+\gamma_{2,1}-k\Lambda\theta_{2,0}-k\Lambda\theta_{2,1})}  \geq 0 \mbox{, }\\ &\gamma_{2,1}+\gamma_{2,0}-k\Lambda\theta_{2,0}-k\Lambda\theta_{2,1} > 0 \mbox{, }\\ &\theta_{2,1}+\theta_{2,0}=1 \mbox{ and }\\ &\theta_{2,1},\theta_{2,0}\geq 0
\end{split}
\end{align}
\normalsize
on the parameters.
Set $\mu = \frac{\gamma_{2,0}-k\Lambda\theta_{2,0}}{\gamma_{2,0}+\gamma_{2,1}-k\Lambda\theta_{2,0}-k\Lambda\theta_{2,1}}$. 
First note that \eqref{eq:twostage2} is equivalent to
\small
\begin{multline}
\label{eq:twostage5}
u_{n+1}  = \argmin_u E_1(u)+\theta_{2,1} L_2(u,U_1)+\theta_{2,0} L_2(u,u_n)+\theta_{2,0}\Lambda \left\| u - u_n \right\|^2+ \\ \theta_{2,1}\Lambda\left\| u - U_1 \right\|^2 + \frac{\gamma_{2,0}+\gamma_{2,1}-k\Lambda\theta_{2,0}-k\Lambda\theta_{2,1}}{2k} \left\| u - \big( \mu u_n + (1-\mu) U_1 \big) \right\|^2.
\end{multline}
\normalsize
This can be seen by expanding the norm squared and comparing
the quadratic and linear terms in $u$. With these tools in hand we can prove energy dissipation: 
\small
\begin{align*}
&E(u_{n+1})\\=&E_1(u_{n+1})+E_2(u_{n+1})\\
\leq& E_1(u_{n+1})+ \theta_{2,1} [L_2(u_{n+1},U_1) +\frac{\Lambda}{2}\left\|u_{n+1}-U_1\right\|^2 ]\\&+ \theta_{2,0} [L_2(u_{n+1},u_n) +\frac{\Lambda}{2}\left\|u_{n+1}-u_n\right\|^2]  && \text{(by \eqref{eq:twostage3})}\\
\leq& E_1(u_{n+1})+ \theta_{2,1} [L_2(u_{n+1},U_1) +\frac{\Lambda}{2}\left\|u_{n+1}-U_1\right\|^2 ]\\&+ \theta_{2,0} [L_2(u_{n+1},u_n) +\frac{\Lambda}{2}\|u_{n+1}-u_n\|^2]\\&+\frac{\gamma_{2,0}+\gamma_{2,1}-k\Lambda\theta_{2,0}-k\Lambda\theta_{2,1}}{2k}\left\| u_{n+1} - \big( \mu u_n + (1-\mu) U_1 \big) \right\|^2.  && \text{(by \eqref{eq:twostage4})}\\
\leq& E_1(U_1) + \theta_{2,1} E_2(U_1)+ \theta_{2,0} [L_2(U_1,u_n) +\frac{\Lambda}{2}\|U_1-u_n\|^2]\\&+\frac{\gamma_{2,0}+\gamma_{2,1}-\Lambda\theta_{2,0}-\Lambda\theta_{2,1}}{2k} \left\|U_1 - \big( \mu u_n + (1-\mu) U_1 \big) \right\|^2 && \text{(by \eqref{eq:twostage5})}\\
\leq& E_1(U_1) + \theta_{2,1}[ L_2(U_1,u_n) +\frac{\Lambda}{2}\|u_{n+1}-u_n\|^2]\\&+ \theta_{2,0} [L_2(U_1,u_n) +\frac{\Lambda}{2}\|U_1-u_n\|^2]\\&+\frac{(\gamma_{2,0}-k\theta_{2,0})^2}{(\gamma_{2,0}+\gamma_{2,1}-k\Lambda\theta_{2,0}-k\Lambda\theta_{2,1})2k} \left\|U_1 - u_n  \right\|^2\\
\leq&   E_1(U_1) + [L_2(U_1,u_n) +\frac{\Lambda}{2}\|U_1-u_n\|^2]+ \frac{\gamma_{1,0}-k\Lambda}{2k} \left\| U_1 - u_n \right\|^2 && \text{(by \eqref{eq:twostage3})}\\
\leq& E(u_n). && \text{(by \eqref{eq:twostage1})}
\end{align*}
\normalsize
The first two conditions of \cref{eq:twostage4} require $k$ to be below a certain threshold. Hence the dissipation of \cref{eq:twostage1} \& \cref{eq:twostage2} is conditional, unless $E_2$ happens to be concave, in which case these two conditions are satisfied for all $k>0$.
\\

We will now extend this discussion to general, $M$-stage case of scheme \cref{eq:ms}:

\begin{theorem}
\label{claim:ms}
Fix a time step $k$. Define $\Lambda=\max\{0,\max_{x,\left\|v\right\|=1} D^2E_2(x)\big(v ,v\big)\}$ and the following auxiliary quantities in terms of the coefficients $\gamma_{m,i}$ and $\theta_{m,i}$ of scheme \cref{eq:ms}:
\begin{align}
\label{eq:tildegamma}
&\tilde{\gamma}_{m,i}=\gamma_{m,i}-k\Lambda\theta_{m,i}-\sum_{j=m+1}^M\tilde{\gamma}_{j,i}\frac{\tilde{S}_{j,m}}{\tilde{S}_{j,j}}\\
\label{eq:S}
&\tilde{S}_{j,m}=\sum_{i=0}^{m-1} \tilde{\gamma}_{j,i}
\end{align}
If $\tilde{S}_{m,m}>0$ for $m=1,\ldots,M$, $\theta_{m-1,i}\geq\theta_{m,i}\geq 0$ and $\sum_{i=0}^{m-1} \theta_{m,i}=1$, then scheme \cref{eq:ms} satisfies the energy stability condition \cref{eq:energymin}: For every $n=0,1,2,\ldots$ we have $E(u_{n+1}) \leq E(u_n)$.
\end{theorem}
As we will see in \cref{sec:ex}, the conditions on the parameters $\gamma_{i,j}$ and $\theta_{m,i}$ of scheme \cref{eq:ms} imposed in \cref{claim:ms} are loose enough to enable meeting consistency conditions to high order.
We will establish \cref{claim:ms} with the help of a couple of lemmas:

\begin{lemma}
Let the auxiliary quantities $\tilde{S}_{j,m}$, and $\tilde{\gamma}_{m,i}$ be defined as in \cref{claim:ms}. 
We have
\label{lem:equiv}
\small
\begin{align*}
&\argmin E(u)+\sum^{m-1}_{i=0} \theta_{m,i} L_2(u,U_i)+\sum^{m-1}_{i=0}\frac{\gamma_{m,i}}{2k}\left\|u-U_i\right\|^2\\=
&\argmin E(u)+\sum^{m-1}_{i=0} \theta_{m,i} [L_2(u,U_i)+\frac{\Lambda}{2}\left\|u-U_i\right\|^2]+\frac{1}{2k}\sum_{j=m}^M \frac{\tilde{S}_{j,m}^2}{\tilde{S}_{j,j}} \left\|u-\sum_{i=0}^{m-1} \frac{\tilde{\gamma}_{j,i}}{\tilde{S}_{j,m}} U_i\right\|^2
\end{align*}
\end{lemma}
\normalsize

\begin{proof}
As in the two step case the proof consists of expanding the norm squared terms and showing that all the quadratic and linear terms of $u$ are equal. First, the expansion of $\sum^{m-1}_{i=0}\frac{\gamma_{m,i}}{2k}\|u-U_i\|^2$ is 
\small
\begin{align}
\label{eq:firstexp}
\frac{\|u\|^2}{2k}\sum^{m-1}_{i=0}\gamma_{m,i} - \frac{1}{k}\langle u,\sum^{m-1}_{i=0}\gamma_{m,i}U_i\rangle+\text{terms that do not depend on $u$.}
\end{align}
\normalsize
Next, we will establish two identities to help us expand
\small\[\frac{1}{2k}\sum_{j=m}^M \frac{\tilde{S}_{j,m}^2}{\tilde{S}_{j,j}} \|u-\sum_{i=0}^{m-1} \frac{\tilde{\gamma}_{j,i}}{\tilde{S}_{j,m}} U_i\|^2.\]\normalsize
First by rearranging \cref{eq:tildegamma},
\small
\begin{equation}
\label{eq:gammaind}
 \gamma_{m,i}-k\Lambda\theta_{m,i}=\sum_{j=m}^M\tilde{\gamma}_{j,i}\frac{\tilde{S}_{j,m}}{\tilde{S}_{j,j}}.
 \end{equation}
 \normalsize
 Next, an identity of $\tilde{S}_{m,m}$:
\small
\begin{align*}
\tilde{S}_{m,m}&=\sum_{i=0}^{m-1} \tilde{\gamma}_{m,i}=\sum_{i=0}^{m-1}\bigg[ \gamma_{m,i}-k\Lambda \theta_{m,i}-\sum_{j=m+1}^M\tilde{\gamma}_{j,i}\frac{\tilde{S}_{j,m}}{\tilde{S}_{j,j}}\bigg]\\&= \sum_{i=0}^{m-1} \bigg[\gamma_{m,i}-k\Lambda \theta_{m,i}\bigg]-\sum_{j=m+1}^M\bigg[\sum_{i=0}^{m-1}\tilde{\gamma}_{j,i}\bigg]\frac{\tilde{S}_{j,m}}{\tilde{S}_{j,j}}\\&=\sum_{i=0}^{m-1}\bigg[ \gamma_{m,i}-k\Lambda \theta_{m,i}\bigg]-\sum_{j=m+1}^M \frac{\tilde{S}_{j,m}^2}{\tilde{S}_{j,j}}.
\end{align*}
\normalsize
We use this identity to establish the following:
\small
\begin{equation}
\label{eq:idenitysum}
    \sum_{j=m}^M \frac{\tilde{S}_{j,m}^2}{\tilde{S}_{j,j}}=\tilde{S}_{m,m}+\sum_{j=m+1}^M \frac{\tilde{S}_{j,m}^2}{\tilde{S}_{j,j}}=\sum_{i=0}^{m-1}\bigg[ \gamma_{m,i}-k\Lambda \theta_{m,i}\bigg]
\end{equation}

\normalsize
Now we can calculate the expansion:
\small
\begin{align*}
&\frac{1}{2k}\sum_{j=m}^M \frac{\tilde{S}_{j,m}^2}{\tilde{S}_{j,j}} \|u-\sum_{i=0}^{m-1} \frac{\tilde{\gamma}_{j,i}}{\tilde{S}_{j,m}} U_i\|^2+\sum^{m-1}_{i=0} \theta_{m,i} \Lambda\left\|u-U_i\right\|^2\\
=&\frac{\|u\|^2}{2k} \sum_{j=m}^M \frac{\tilde{S}_{j,m}^2}{\tilde{S}_{j,j}} -\frac{1}{k}\langle u,\sum_{i=0}^{m-1} \sum_{j=m}^M \tilde{\gamma}_{j,i} \frac{\tilde{S}_{j,m}}{\tilde{S}_{j,j}} U_i \rangle+\frac{\Lambda}{2}\|u\|^2\sum^{m-1}_{i=0} \theta_{m,i} +\Lambda \sum^{m-1}_{i=0} \langle u,\theta_{m,i}U_i\rangle\\&+\text{terms that do not depend on $u$}\\
=&\frac{\|u\|^2}{2k}\sum^{m-1}_{i=0}\gamma_{m,i} - \frac{1}{k}\langle u,\sum^{m-1}_{i=0}\gamma_{m,i}U_i\rangle+\text{terms that do not depend on $u$.}
\end{align*}\\
\normalsize Where the last equality follows from \cref{eq:gammaind} and \cref{eq:idenitysum}. Since this expansion matches \cref{eq:firstexp} up to a constant in $u$ the proof is complete.
\end{proof}

\begin{lemma}
\label{lem:step}
Let $\Lambda$ and the auxiliary quantities $\tilde{S}_{j,m}$, $\tilde{\gamma}_{m,i}$ be given in \cref{claim:ms}. Additionally, let $\tilde{S}_{m,m}>0$ for $m=1,\ldots,M$. Then
\small
\begin{align*}
    &E_1(U_m)+\sum^{m-1}_{i=0} \theta_{m,i} [L_2(U_m,U_i)+\frac{\Lambda}{2}\left\|U_m-U_i\right\|^2]+\frac{1}{2k}\sum_{j=m}^M \frac{\tilde{S}_{j,m}^2}{\tilde{S}_{j,j}} \left\|U_m-\sum_{i=0}^{m-1} \frac{\tilde{\gamma}_{j,i}}{\tilde{S}_{j,m}} U_i\right\|^2\\\leq & E_1(U_{m-1})+\sum^{m-2}_{i=0} \theta_{m-1,i} [L_2(U_{m-1},U_i)+\frac{\Lambda}{2}\left\|U_{m-1}-U_i\right\|^2]\\&+\frac{1}{2k}\sum_{j=m-1}^M \frac{\tilde{S}_{j,m-1}^2}{\tilde{S}_{j,j}} \left\|U_{m-1}-\sum_{i=0}^{m-2} \frac{\tilde{\gamma}_{j,i}}{\tilde{S}_{j,m-1}} U_i\right\|^2
\end{align*}
\normalsize
\end{lemma}

\begin{proof}

By \cref{eq:ms} \& \cref{lem:equiv},
\small
\[
U_m=\argmin_u E(u)+\sum^{m-1}_{i=0} \theta_{m,i} [L_2(u,U_i)+\Lambda\left\|u-U_i\right\|^2]+\frac{1}{2k}\sum_{j=m}^M \frac{\tilde{S}_{j,m}^2}{\tilde{S}_{j,j}} \|u-\sum_{i=0}^{m-1} \frac{\tilde{\gamma}_{j,i}}{\tilde{S}_{j,m}} U_i\|^2.
\]
\normalsize
Since $U_m$ is the minimizer of the above optimization problem
\small
\begin{align}
&E_1(U_m)+\sum^{m-1}_{i=0} \theta_{m,i} [L_2(U_m,U_i)+\frac{\Lambda}{2}\left\|U_m-U_i\right\|^2]\nonumber\\&+\frac{1}{2k}\sum_{j=m}^M \frac{\tilde{S}_{j,m}^2}{\tilde{S}_{j,j}} \left\|U_m-\sum_{i=0}^{m-1} \frac{\tilde{\gamma}_{j,i}}{\tilde{S}_{j,m}} U_i\right\|^2 \nonumber\\
\label{eq:lem2step1}
\begin{split}
\leq&E_1(U_{m-1})+\theta_{m,m-1}E_2(U_{m-1})+\sum^{m-2}_{i=0} \theta_{m,i} [L_2(U_{m-1},U_i)+\frac{\Lambda}{2}\left\|U_{m-1}-U_i\right\|^2]\\&+\frac{1}{2k}\sum_{j=m}^M \frac{\tilde{S}_{j,m}^2}{\tilde{S}_{j,j}} \left\|U_{m-1}-\sum_{i=0}^{m-1} \frac{\tilde{\gamma}_{j,i}}{\tilde{S}_{j,m}} U_i\right\|^2
\end{split}
\end{align}
\normalsize
We give two inequalities to aid us in the proof. 
First, using the definition of the auxiliary variables, we can state an identity that will simplify \cref{eq:lem2step1}. For $m>1$ and $j\geq m$ 
\small
\begin{align}
\label{eq:lem2step12}
\begin{split}
&\frac{\tilde{S}_{j,m}^2}{\tilde{S}_{j,j}} \left\|U_{m-1}-\sum_{i=0}^{m-1} \frac{\tilde{\gamma}_{j,i}}{\tilde{S}_{j,m}} U_i\right\|^2=\frac{\tilde{S}_{j,m}^2}{\tilde{S}_{j,j}} \left\|U_{m-1}\bigg(1-\frac{\tilde{\gamma}_{j,m-1}}{\tilde{S}_{j,m}}\bigg)-\sum_{i=0}^{m-2} \frac{\tilde{\gamma}_{j,i}}{\tilde{S}_{j,m}} U_i\right\|^2\\
=&\frac{\tilde{S}_{j,m}^2}{\tilde{S}_{j,j}} \left\|U_{m-1}\bigg(\frac{\tilde{S}_{j,m-1}}{\tilde{S}_{j,m}}\bigg)-\sum_{i=0}^{m-2} \frac{\tilde{\gamma}_{j,i}}{\tilde{S}_{j,m}} U_i\right\|^2=\frac{\tilde{S}_{j,m-1}^2}{\tilde{S}_{j,j}} \left\|U_{m-1}-\sum_{i=0}^{m-2} \frac{\tilde{\gamma}_{j,i}}{\tilde{S}_{j,m-1}} U_i\right\|^2.
\end{split}
\end{align}
\normalsize
Now since $\tilde{S}_{m-1,m-1}>0$,
\small
\begin{equation}
   \label{eq:lem2step3} 
    \frac{\tilde{S}_{m-1,m-1}^2}{\tilde{S}_{m-1,m-1}} \left\|U_{m-1}-\sum_{i=0}^{m-2} \frac{\tilde{\gamma}_{m-1,i}}{\tilde{S}_{m-1,m-1}} U_i\right\|^2>0.
\end{equation}
\normalsize
Using \cref{eq:lem2step12} and \cref{eq:lem2step3} we have
\small
\begin{multline}
\label{eq:lem2step2}
\frac{1}{2k}\sum_{j=m}^M \frac{\tilde{S}_{j,m}^2}{\tilde{S}_{j,j}} \left\|U_{m-1}-\sum_{i=0}^{m-1} \frac{\tilde{\gamma}_{j,i}}{\tilde{S}_{j,m}} U_i\right\|^2=\frac{1}{2k}\sum_{j=m}^M \frac{\tilde{S}_{j,m-1}^2}{\tilde{S}_{j,j}} \left\|U_{m-1}-\sum_{i=0}^{m-2} \frac{\tilde{\gamma}_{j,i}}{\tilde{S}_{j,m-1}} U_i\right\|^2\\\leq\frac{1}{2k}\sum_{j=m-1}^M \frac{\tilde{S}_{j,m-1}^2}{\tilde{S}_{j,j}} \left\|U_{m-1}-\sum_{i=0}^{m-2} \frac{\tilde{\gamma}_{j,i}}{\tilde{S}_{j,m-1}} U_i\right\|^2.
\end{multline}
\normalsize
Next, since $\sum_{i=1}^{m-1}\theta_{m,i}=1$ for all $m$ we have the equality 
\begin{equation}
\label{eq:thetaeq}
   \theta_{m,m-1}=1-\sum^{m-2}_{i=0} \theta_{m,i}=\sum^{m-2}_{i=0} \theta_{m-1,i}-\sum^{m-2}_{i=0}\theta_{m,i} 
\end{equation} 
Using \cref{eq:twostage3} and \cref{eq:thetaeq}, we have our second inequality:
\small
\begin{align}
\label{eq:lem2stepE}
\begin{split}
&\theta_{m,m-1}E_2(U_{m-1})+\sum^{m-2}_{i=0} \theta_{m,i} [L_2(U_{m-1},U_i)+\frac{\Lambda}{2}\left\|U_{m-1}-U_i\right\|^2]\\
=&\sum^{m-2}_{i=0} (\theta_{m-1,i}-\theta_{m,i})E_2(U_{m-1}) +\sum^{m-2}_{i=0} \theta_{m,i} [L_2(U_{m-1},U_i)+\frac{\Lambda}{2}\left\|U_{m-1}-U_i\right\|^2]\\
\leq&\sum^{m-2}_{i=0} (\theta_{m-1,i}-\theta_{m,i}) [L_2(U_{m-1},U_i)+\frac{\Lambda}{2}\left\|U_{m-1}-U_i\right\|^2]\\&+\sum^{m-2}_{i=0} \theta_{m,i} [L_2(U_{m-1},U_i)+\frac{\Lambda}{2}\left\|U_{m-1}-U_i\right\|^2]\\
=&\sum^{m-2}_{i=0} \theta_{m-1,i} [L_2(U_{m-1},U_i)+\frac{\Lambda}{2}\left\|U_{m-1}-U_i\right\|^2]
\end{split}
\end{align}
\normalsize
Using inequalities \cref{eq:lem2step2} and \cref{eq:lem2stepE}, we have that \cref{eq:lem2step1} is less than or equal to
\small
\begin{multline*}
E_1(U_{m-1})+\sum^{m-2}_{i=0} \theta_{m-1,i} [L_2(U_{m-1},U_i)+\frac{\Lambda}{2}\left\|U_{m-1}-U_i\right\|^2]\\+\frac{1}{2k}\sum_{j=m-1}^M \frac{\tilde{S}_{j,m-1}^2}{\tilde{S}_{j,j}} \left\|U_{m-1}-\sum_{i=0}^{m-2} \frac{\tilde{\gamma}_{j,i}}{\tilde{S}_{j,m-1}} U_i\right\|^2
\end{multline*}
\normalsize
concluding the proof.
\end{proof}

\begin{proof}(of theorem)
The main idea of the proof is to use \cref{lem:step} repeatedly to relate the energy of $E(u_{n+1})$ to $E(u_n)$. First, by \cref{eq:twostage3} and our assumption that $\tilde{S}_{M,M}>0$
\small
\begin{align*}
E(u_{n+1})&=E_1(U_M)+E_2(U_M)\\
&\leq E_1(U_M)+\sum^{M-1}_{i=0} \theta_{M,i} [L_2(U_M,U_i)+\frac{\Lambda}{2}\left\|U_M-U_i\right\|^2]\\&+\frac{1}{2k} \frac{\tilde{S}_{M,M}^2}{\tilde{S}_{M,M}} \|U_M-\sum_{i=0}^{M-1} \frac{\tilde{\gamma}_{M,i}}{\tilde{S}_{M,M}} U_i\|^2.
\end{align*}
\normalsize
By using the \cref{lem:step} repeatedly we have
\small
\begin{align*}
&E_1(U_M)+\sum^{M-1}_{i=0} \theta_{M,i} [L_2(U_M,U_i)+\frac{\Lambda}{2}\left\|U_M-U_i\right\|^2]+\frac{1}{2k} \frac{\tilde{S}_{M,M}^2}{\tilde{S}_{M,M}} \|U_M-\sum_{i=0}^{M-1} \frac{\tilde{\gamma}_{M,i}}{\tilde{S}_{M,M}} U_i\|^2\\
\leq& E_1(U_{M-1})+\sum^{M-2}_{i=0} \theta_{M-1,i} [L_2(U_{M-1},U_i)+\frac{\Lambda}{2}\left\|U_{M-1}-U_i\right\|^2]\\&+\frac{1}{2k}\sum_{j=M-1}^M \frac{\tilde{S}_{j,M-1}^2}{\tilde{S}_{j,j}} \left\|U_{M-1}-\sum_{i=0}^{M-2} \frac{\tilde{\gamma}_{j,i}}{\tilde{S}_{j,M-1}} U_i\right\|^2\\
& \vdots \\\leq& E_1(U_1)+ L_2(U_1,U_0)+\frac{\Lambda}{2}\left\|U_1-U_0\right\|^2+\frac{1}{2k}\sum_{j=1}^M \frac{\tilde{S}_{j,1}^2}{\tilde{S}_{j,j}} \|U_1- \frac{\tilde{\gamma}_{j,0}}{\tilde{S}_{j,1}} U_0\|^2.
\end{align*}
\normalsize 
By \cref{eq:ms} and \cref{lem:equiv}
\small
\[U_1=\argmin_u E_1(u)+L_2(u,U_0)+\frac{\Lambda}{2}\left\|u-U_0\right\|^2+\frac{1}{2k}\sum_{j=1}^M \frac{\tilde{S}_{j,1}^2}{\tilde{S}_{j,j}} \|u- \frac{\tilde{\gamma}_{j,0}}{\tilde{S}_{j,1}} U_0\|^2\]
\normalsize
so
\small
\begin{align*}
&E_1(U_1)+L_2(U_1,U_0)+\frac{\Lambda}{2}\left\|U_1-U_0\right\|^2+\frac{1}{2k}\sum_{j=1}^M \frac{\tilde{S}_{j,1}^2}{\tilde{S}_{j,j}} \|U_1- \frac{\tilde{\gamma}_{j,0}}{\tilde{S}_{j,1}} U_0\|^2\\ \leq& E_1(U_0)+E_2(U_0)+\frac{\Lambda}{2}\left\|U_0-U_0\right\|^2+\frac{1}{2k}\sum_{j=1}^M \frac{\tilde{S}_{j,1}^2}{\tilde{S}_{j,j}} \|U_0-U_0\|^2\\=&E(u_n)
\end{align*}
\normalsize
completing the proof of the theorem.
\end{proof}

\begin{remark}
In the above proof, we assume that $E_2(u)$ is two times differentiable.
This assumption can be dropped if we replace $L_2(u,p)$ with another approximation $A_2(u,p)$ that has the properties $A_2(u,u)=E_2(u)$ and for some choice $\Lambda$, $E_2(u) \leq A_2(u,p)+\frac{\Lambda}{2}\left\|u-p\right\|^2$ for all $u$ and $p$.
\end{remark}
\section{Examples of the New Schemes for Gradient Flows}
\label{sec:ex}
In this section, we give examples of {\it high order} semi-implicit schemes for gradient flows, for any desired choice of implicit and explicit terms $E_1$ and $E_2$, that are energy stable under the conditions of \cref{claim:ms}.
First, we give the conditions on $\gamma_{m,i}$ and $\theta_{m,i}$ in scheme \cref{eq:ms} to ensure high order consistency with the abstract evolution law \cref{eq:de}.
Recall that $U_0=u_n$. From \cref{eq:ms}, each stage $U_m$ satisfies the Euler-Lagrange equation:
\begin{equation}
    \label{eq:ms2}
    \bigg[\sum^{m-1}_{i=0}\gamma_{m,i}\bigg]U_{m}+k\nabla_HE_1(U_m)=-\sum^{m-1}_{i=0} k\theta_{m,i}\nabla_HE_2(U_i)+\sum^{m-1}_{i=0}\gamma_{m,i}U_i.
\end{equation}
\cref{eq:ms2} is equivalent to the form more often seen for ARK IMEX methods:
\small
\begin{equation}
\label{eq:rk}
U_m=U_0-k\sum^{m}_{i=1} \alpha_{m,i} \nabla_H E_1(U_i)-k\sum^{m-1}_{i=1} \tilde{\alpha}_{m,i} \nabla_H E_2(U_i)
\end{equation}
\normalsize
where $\alpha_{m,i}$ and $\tilde{\alpha}_{m,i}$ depend on $\gamma_{m_i}$ and $\theta_{m,i}$.
The consistency equations for ARK IMEX methods have been previously worked out \cite{kennedy2003additive,Pareschi2005,shin2017unconditionally,zharovsky2015class}.
As such, we will state without proof the conditions required to achieve various orders of accuracy in terms of $\gamma$ and $\theta$:

\begin{claim}
\label{claim:cons}
Let $U_i$ be given in \cref{eq:ms}. The Taylor expansion of $U_i$ at each stage has the form: 
\small
\begin{align}
\label{eq:indte}
\begin{split}
    U_i&=U_0-\beta_{1,i}k DE(U_0)+k^2\big[\beta_{2,i}k^2D^2E_1(U_0)DE(U_0)+\beta_{3,i}D^2E_2(U_0)DE(U_0)\big]\\&-k^3\big[\beta_{4,i}D^2E_1(U_0)\left(D^2E_1(U_0)\left(DE(U_0)\right)\right)+\beta_{5,i}D^2E_1(U_0)\left(D^2E_2(U_0)\left(DE(U_0)\right)\right)\\&+\beta_{6,i}D^2E_2(U_0)\left(D^2E_1(U_0)\left(DE(U_0)\right)\right)+\beta_{7,i}D^2E_2(U_0)\left(D^2E_2(U_0)\left(DE(U_0)\right)\right)\\&+\beta_{8,i}D^3E_1(U_0)\big(DE(U_0),DE(U_0)\big)+\beta_{9,i}D^3E_2(U_0)\big(DE(U_0),DE(U_0)\big)\big]+\text{h.o.t.}
    \end{split}
\end{align}
\normalsize
where for $l\in\{1,2,3,\ldots\}$, $D^l E(u) : H^l \to \mathbb{R}$ denotes the multilinear form given by
\small
\[ D^l E(u)\big(v_1,\ldots,v_n\big) = \left. \frac{\partial^l}{\partial s_1 \cdots \partial s_l} E(u+s_1 v_1 + s_2 v_2 + \cdots + s_l v_l) \right|_{s_1 = s_2 = \cdots = s_l = 0}\] \normalsize 
so that the linear functional $D^lE(u)\big(v_1,v_2,\ldots,v_{l-1},\cdot \big) : H \to \mathbb{R}$ may be identified with an element of $H$, and so on.
The coefficients of \cref{eq:indte} obey the following recursive relations:
\begin{align}
\begin{split}
\label{eq:rec}
&\beta_{1,0}=\beta_{2,0}=\ldots=\beta_{9,0}=0\\
&\beta_{1,m}=\frac{1}{S_m}\bigg[1+\sum^{m-1}_{i=1}\gamma_{m,i}\beta_{1,i} \bigg]\\
&\beta_{2,m}=\frac{1}{S_m}\bigg[\beta_{1,m}+\sum^{m-1}_{i=1}\gamma_{m,i}\beta_{2,i} \bigg]\\
&\beta_{3,m}=\frac{1}{S_m}\bigg[\sum_{i=0}^{m-1} \theta_{m,i} \beta_{1,i}+\sum^{m-1}_{i=1}\gamma_{m,i}\beta_{3,i} \bigg]\\
&\beta_{4,m}=\frac{1}{S_m}\bigg[\beta_{2,m}+\sum^{m-1}_{i=1}\gamma_{m,i} \beta_{4,i}\bigg] \\
&\beta_{5,m}=\frac{1}{S_m}\bigg[\beta_{3,m}+\sum^{m-1}_{i=1}\gamma_{m,i} \beta_{5,i}\bigg] \\
&\beta_{6,m}=\frac{1}{S_m}\bigg[\sum_{i=0}^{m-1} \theta_{m,i} \beta_{2,i}+\sum^{m-1}_{i=1}\gamma_{m,i}\beta_{6,i} \bigg]\\
&\beta_{7,m}=\frac{1}{S_m}\bigg[\sum_{i=0}^{m-1} \theta_{m,i} \beta_{3,i}+\sum^{m-1}_{i=1}\gamma_{m,i}\beta_{7,i} \bigg]\\
&\beta_{8,m}=\frac{1}{S_m}\bigg[\frac{\beta_{1,m}^2}{2}+\sum^{m-1}_{i=1}\gamma_{m,i} \beta_{8,i}\bigg] \\
&\beta_{9,m}=\frac{1}{S_m}\bigg[\frac{1}{2}\sum_{i=0}^{m-1} \theta_{m,i} \beta_{1,i}^2+\sum^{m-1}_{i=1}\gamma_{m,i}\beta_{9,i} \bigg]
\end{split}
\end{align}
with $S_m=\sum^{m-1}_{i=0}\gamma_{m,i}$.
Furthermore, the following conditions for $u_{n+1}=U_M$ in scheme \cref{eq:ms} are necessary and
sufficient for various orders of accuracy:

\begin{alignat}{5}
\label{eq:cons}
 &\text{\underline{First Order:}}&\quad & \text{\underline{Second Order:}} &\quad &\text{\underline{Third Order:}} \nonumber \\
&\beta_{1,M}=1& &\beta_{1,M}=1& &\beta_{1,M}=1 \nonumber\\
& & &\beta_{2,M}=1/2 & &\beta_{2,M}=1/2\\
& & &\beta_{3,M}=1/2& &\beta_{3,M}=1/2 \nonumber\\
& & & & &\beta_{4,M}=\beta_{5,M}=\ldots=\beta_{9,M}=1/6 \nonumber
\end{alignat}
\end{claim}

Now, we give second order and a third order example of method \cref{eq:ms}. However, The examples we give are not unique by any means. We begin with a five step method that is second order accurate:

 \begin{align}
\label{eq:2ndordergamma}
\begin{split}
 &\theta \approx \left(
\begin{array}{ccccc}
 1. & 0 & 0 & 0 & 0 \\
 0.009 & 0.991 & 0 & 0 & 0 \\
 0.009 & 0.991 & 0 & 0 & 0 \\
 0 & 0 & 0 & 1. & 0 \\
 0 & 0 & 0 & 1. & 0 \\
\end{array}
\right)\\
&\gamma \approx \left(\begin{array}{ccccc} 8.841 & 0 & 0 & 0 & 0\\ -0.925 & 5.360 & 0 & 0 & 0\\ -4.443 & 6.041 & 0.950 & 0 & 0\\ -3.288 & 5.895 & -0.351 & 0.172 & 0\\ -3.895 & -0.335 & 4.964 & -1.722 & 7.684 \end{array}\right)
\end{split}
 \end{align}
 
which is stable for $k\Lambda \leq 3/872$.

Next we have a thirteen step method that is third order accurate:
\tiny
\begin{align}
\label{eq:3rdordergamma}
\begin{split}
 &\theta\approx \left(
\begin{array}{ccccccccccccc}
 1. & 0 & 0 & 0 & 0 & 0 & 0 & 0 & 0 & 0 & 0 & 0 & 0 \\
 0.049 & 0.951 & 0 & 0 & 0 & 0 & 0 & 0 & 0 & 0 & 0 & 0 & 0 \\
 0.024 & 0.075 & 0.901 & 0 & 0 & 0 & 0 & 0 & 0 & 0 & 0 & 0 & 0 \\
 0.017 & 0.042 & 0.113 & 0.829 & 0 & 0 & 0 & 0 & 0 & 0 & 0 & 0 & 0 \\
 0.012 & 0.029 & 0.071 & 0.386 & 0.501 & 0 & 0 & 0 & 0 & 0 & 0 & 0 & 0 \\
 0.01 & 0.023 & 0.06 & 0.366 & 0.457 & 0.085 & 0 & 0 & 0 & 0 & 0 & 0 & 0 \\
 0.007 & 0.018 & 0.05 & 0.351 & 0.437 & 0.06 & 0.076 & 0 & 0 & 0 & 0 & 0 & 0 \\
 0.003 & 0.005 & 0.006 & 0.008 & 0.009 & 0.011 & 0.028 & 0.929 & 0 & 0 & 0 & 0 & 0 \\
 0.002 & 0.002 & 0.002 & 0.002 & 0.003 & 0.004 & 0.009 & 0.029 & 0.948 & 0 & 0 & 0 & 0 \\
 0 & 0.001 & 0.001 & 0.001 & 0.001 & 0.002 & 0.004 & 0.007 & 0.011 & 0.971 & 0 & 0 & 0 \\
 0 & 0 & 0.001 & 0.001 & 0.001 & 0.001 & 0.003 & 0.005 & 0.008 & 0.912 & 0.069 & 0 & 0 \\
 0 & 0 & 0 & 0 & 0 & 0.001 & 0.002 & 0.003 & 0.005 & 0.107 & 0.025 & 0.857 & 0 \\
 0 & 0 & 0 & 0 & 0 & 0 & 0.001 & 0.001 & 0.002 & 0.013 & 0.007 & 0.018 & 0.958 \\
\end{array}
\right)\\
&\gamma \approx \left(
\begin{array}{ccccccccccccc}
 11. & 0 & 0 & 0 & 0 & 0 & 0 & 0 & 0 & 0 & 0 & 0 & 0 \\
 2.1 & 15.5 & 0 & 0 & 0 & 0 & 0 & 0 & 0 & 0 & 0 & 0 & 0 \\
 1.4 & 1.6 & 17. & 0 & 0 & 0 & 0 & 0 & 0 & 0 & 0 & 0 & 0 \\
 0.2 & 1.6 & -2.4 & 18.1 & 0 & 0 & 0 & 0 & 0 & 0 & 0 & 0 & 0 \\
 0.3 & -8.5 & 3. & 9.6 & 7.8 & 0 & 0 & 0 & 0 & 0 & 0 & 0 & 0 \\
 -1.4 & -5.9 & -0.1 & 2. & 8. & 4.1 & 0 & 0 & 0 & 0 & 0 & 0 & 0 \\
 -4. & -0.5 & -0.4 & -1.8 & 5.1 & 6.8 & 0.9 & 0 & 0 & 0 & 0 & 0 & 0 \\
 -9.2 & 4.8 & 2.7 & -3.2 & 2.5 & 6.2 & 2.5 & 4.6 & 0 & 0 & 0 & 0 & 0 \\
 -1.7 & -3.6 & -0.1 & 1.3 & 5.7 & 3.4 & -0.8 & -0.8 & 0.4 & 0 & 0 & 0 & 0 \\
 -2.7 & -3.5 & 0.6 & 1.4 & 6.1 & 3.5 & -0.7 & -0.2 & -0.4 & 0.5 & 0 & 0 & 0 \\
 5.9 & -4.8 & -5.1 & -3.1 & 3.4 & 6.6 & -0.7 & -5.2 & 4.9 & -0.8 & 8.2 & 0 & 0 \\
 7.1 & 0.9 & -3.1 & -2.7 & -5.8 & -1.9 & 0.6 & -3.4 & 4.3 & -1.3 & 9.2 & 9.1 & 0 \\
 3.8 & 1.9 & 2.7 & 2.1 & -7.5 & -10.6 & -1.2 & 2. & 0.7 & -0.2 & -0.2 & 9.5 & 12.8 \\
\end{array}
\right)
\end{split}
 \end{align}
\normalsize
and is stable if $k\Lambda \leq 18/28567$. The coefficients to machine precision as well as code to verify \cref{claim:ms} and \cref{claim:cons} can be found at \url{https://github.com/AZaitzeff/SIgradflow}.
In the following section, we consider methods for \cref{eq:Lde}, when the inner product changes with the solution.

\section{Schemes for Solving Gradient Flows with Solution Dependent Inner Product}
\label{sec:SDIP}
Now we move on to the problem of simulating flow  
\cref{eq:Lde},
\begin{equation*}
    u'=-\op(u)\nabla_H E(u).
\end{equation*}
We consider the case where $\op(u)$ is strictly positive definite. Our approach will be as follows: 
\begin{enumerate}
    \item Generate a $u_*$ from $u_n$.
    \item Construct $\op(u_*)$.
    \item Use the algorithm \cref{eq:ms} with norm $\left\|\cdot\right\|^2_{\op^{-1}(u_*)}= \langle \cdot ,\op^{-1}(u_*) \cdot \rangle$ to generate $u_{n+1}$.
\end{enumerate}

One advantage to constructing $\op(u_*)$ and then using it in \cref{eq:ms} is that \cref{claim:ms} immediately gives conditional energy stability for coefficients such as \cref{eq:2ndordergamma} or \cref{eq:3rdordergamma}. Thus, we only need to consider what choice of $u_*$ will give our algorithm the desired level of accuracy.
Now at every step we are solving 
\small
\begin{equation}
    \label{eq:Lms2}
    \bigg[\sum^{m-1}_{i=0}\gamma_{m,i}\bigg]U_{m}+k\op(u_*)\nabla_HE_1(U_m)=-k\op(u_*)\sum^{m-1}_{i=0} \theta_{m,i}\nabla_HE_2(U_i)+\sum^{m-1}_{i=0}\gamma_{m,i}U_i.
\end{equation}
\normalsize
We will set up the consistency equations for \cref{eq:Lde}.
Let $u_n=u(t_0)$. For convenience, denote $\op(u_n)$ as $\op_n$ and $E(u_n)$ as $E_n$. We begin with the exact solution starting from $u(t_0)$:
\small
\[
\begin{cases}
u_t=-\op(u)\nabla E(u) & t>t_0 \\
u(t_0)=U_0
  \end{cases}
\]
\normalsize
By Taylor expanding around $t_0$ we find
\begin{equation}
\label{eq:truewmu}
u(k+t_0)=u(t_0)+ku_t(t_0)+\frac{1}{2}k^2u_{tt}(t_0)+\frac{1}{6}k^3u_{ttt}(t_0)
\end{equation}
where the higher derivatives in time are found using \cref{eq:Lde}:
\small
\begin{align*}
    u_t(t_0)=&- \op_nDE_n\\
    u_{tt}(t_0)=& D\op_n(\op_nDE_n) DE_n+\op_n D^2E_n(\op_nDE_n)\\
    u_{ttt}(t_0)=&- D\op_n(D\op_n(\op_nDE_n)DE_n) D E_n - D^2\op_n(\op_nDE_n,\op_nDE_n) DE_n\\
    &- D\op_n(\op_n(D^2E_n(\op_nDE_n))) D E_n-2D\op_n(\op_nDE_n)D^2E_n (\op_nDE_n)\\&-\op_nD^2E_n(D\op_n(\op_nDE_n)) DE_n-\op_nD^2E_n(\op_nD^2E_n(\op_nDE_n))\\&-\op_n
    D^3E_n\big(\op_nDE_n,\op_nDE_n\big)
\end{align*}
\normalsize
where for $l\in\{1,2,3,\ldots\}$, $D^l L(u) : H^l \to H$ denotes the multilinear form given by
\small
\[ D^l \op(u)\big(v_1,\ldots,v_l\big) = \left. \frac{\partial^l}{\partial s_1 \cdots \partial s_l} \op(u+s_1 v_1 + s_2 v_2 + \cdots + s_l v_l) \right|_{s_1 = s_2 = \cdots = s_l = 0}\] \normalsize 
so that $D^l\op(u)\big(v_1,v_2,\ldots,v_{l}\big)$ is a linear operator from $H$ to $H$.

In the next two subsections, we provide second and third order examples and accompanying consistency calculations. Both of these examples also have the property that $E(u_*) \leq E(u_n)$.
\subsection{Second Order Method}
\begin{algorithm}[H]
\caption{A second order method for solving gradient flows with solution dependent inner product}
\label{alg:m2ndorder}
Fix a time step size $k>0$.
Set $u_n = u_0$.
To obtain $u_{n+1}$ from $u_n$, carry out the following steps:
\begin{enumerate}
\item Find $u_*$ by solving 
    $u_*+\frac{1}{2}k\op_n\nabla E_1(u_*)=u_n-\frac{1}{2}k\op_n\nabla E_2(u_n)$
\item Find $u_{n+1}$ using \cref{eq:Lms2} with coefficients \cref{eq:2ndordergamma} and $u_*$ in $\op(u_*)$ given by Step 1 of this algorithm.
\end{enumerate}
\end{algorithm}

Our second order algorithm is laid out in \cref{alg:m2ndorder}. Now we will prove that it is indeed second order. First, the expansion of $u_*$ is
\small
\begin{align}
\begin{split}
\label{eq:1s2o}
    u_*&=u_n-\frac{1}{2}k\op_n DE_n+O(k^2)
    \end{split}
\end{align}
\normalsize
We use \cref{eq:rec} to get an expansion of $u_{n+1}$:
\small
\begin{equation}
\label{eq:Lsecondorder}
    u_{n+1}=u_n-k\op(u_*) DE_n+\frac{1}{2}k^2\op(u_*)D^2E_n(\op(u_*)DE_n)+O(k^3)
\end{equation}
\normalsize
Now, expand $u_*$ around $u_n$ in \cref{eq:Lsecondorder}:
\begin{align*}
    u_{n+1}=&u_n-k\op_n DE_n-kD\op_n(u_*-u_n) DE_n\\&+\frac{1}{2}k^2\op_nD^2E_n(\op_nDE_n)+O(k^3)\\
    =&u_n-k\op_n DE_n+\frac{1}{2}k^2D\op_n(\op_n DE_n) DE_n\\&+\frac{1}{2}k^2\op_nD^2E_n(\op_n DE_n)+O(k^3)
\end{align*}
The Taylor expansion of $u_{n+1}$ matches \cref{eq:truewmu} to second order.
\subsection{Third Order Method}

\begin{algorithm}[H]
\caption{A third order method for solving gradient flows with solution dependent inner product}
\label{alg:m3rdorder}
Fix a time step size $k>0$ and
set $u_n = u_0$. For convenience, we will denote
$D^2\op(u_*)\big(\op(u_*)\nabla E(u_*),\op(u_*)\nabla E(u_*)\big)$ as $D^2\op(u_*)$. Additionally, let $MS \Big(\tilde{k},\op(u_*),\tilde{u},\gamma,\theta\Big)$ denote $U_M$ obtained from the multistage algorithm
\small
\[
\bigg[\sum^{m-1}_{i=0}\gamma_{m,i}\bigg]U_{m}+\tilde{k}\op(u_*)\nabla_HE_1(U_m)=-\tilde{k}\op(u_*)\sum^{m-1}_{i=0} \theta_{m,i}\nabla_HE_2(U_i)+\sum^{m-1}_{i=0}\gamma_{m,i}U_i.
\]
\normalsize
with $U_0 = \tilde{u}$.
To obtain $u_{n+1}$ from $u_n$, carry out the following steps:
\begin{enumerate}
\item  Let $\gamma$ and $\theta$ be given by \cref{eq:si1c}. Set $u_{*_1}=MS\Big(\frac{1}{6}k,\op(u_n),u_n,\gamma,\theta\Big).$
\item Let $\gamma$ and $\theta$ be given by \cref{eq:3rdordergamma}. Set \small\[\bar{u}=MS\bigg(\frac{1}{2}k,\op(u_{*_1})-\frac{1}{72}k^2D^2\op(u_{*_1}),u_n,\gamma,\theta\bigg).\]\normalsize
\item Let $\gamma=\theta=\left(1\right)$. Set $u_{*_{2,1}}=MS\Big(\frac{2}{5}k,\op(u_n),u_n,\gamma,\theta\Big)$.
\item Let $\gamma$ and $\theta$ be given by \cref{eq:si1125c}. Set $u_{*_{2,2}}=MS\Big(\frac{5}{6}k,\op(u_{*_{2,1}}),u_n,\gamma,\theta\Big).$
\item Let $\gamma$ and $\theta$ be given by \cref{eq:3rdordergamma}. Then \small\[u_{n+1}=MS\bigg(\frac{1}{2}k,\op(u_{*_{2,2}})-\frac{1}{72}k^2D^2\op(u_{*_{2,2}}),\bar{u},\gamma,\theta\bigg).\]\normalsize
\end{enumerate}
\end{algorithm}
Now we present our third order algorithm for solving \cref{eq:Lde}. It requires the use of two new sets of coefficients, 
\small
\begin{align}
\label{eq:si1c}
\begin{split}
    &\theta \approx\left(
\begin{array}{ccc}
 1. & 0 & 0.  \\
 -0.667 & 0.333 & 0 \\
 0 & 0 & 1.000\\
\end{array}
\right)\\
    &\gamma \approx\left(
\begin{array}{ccc}
 1.833 & 0 & 0.  \\
 0.556 & 0.667 & 0 \\
 1.030 & -0.026 & 0.159 \\
\end{array}
\right)
\end{split}
 \end{align}
\normalsize
and
\small
\begin{align}
\label{eq:si1125c}
\begin{split}
    &\theta \approx\left(
\begin{array}{ccccccc}
 1. & 0 & 0 & 0 & 0 & 0 & 0 \\
 0.708 & 0.292 & 0 & 0 & 0 & 0 & 0 \\
 0.013 & 0.018 & 0.969 & 0 & 0 & 0 & 0 \\
 0.008 & 0.012 & 0.867 & 0.113 & 0 & 0 & 0 \\
 0.006 & 0.009 & 0.206 & 0.056 & 0.724 & 0 & 0 \\
 0 & 0.005 & 0.05 & 0.025 & 0.053 & 0.867 & 0 \\
 0 & 0 & 0.015 & 0.009 & 0.015& 0.04 & 0.920 \\
\end{array}
\right)\\
    &\gamma \approx \left(
\begin{array}{ccccccc}
 7.727 & 0 & 0 & 0 & 0 & 0 & 0 \\
 0.594 & 2.241 & 0 & 0 & 0 & 0 & 0 \\
 3.056 & -0.455 & 0.636 & 0 & 0 & 0 & 0 \\
 -1.571 & 5.091 & -1.063 & 2.786 & 0 & 0 & 0 \\
 -3.714 & 3.1 & -1.267 & 1.545 & 9.655 & 0 & 0 \\
 -6.923 & 5.1 & -2.056 & 3.471 & 4.571 & 4.033 & 0 \\
 -2.467 & -2.1 & 0.009 & -0.182 & 0.660 & 7.224 & 9.428 \\
\end{array}
\right),
\end{split}
 \end{align}
 \normalsize
 to achieve particular Taylor expansions as we explain later in the section. The values of \cref{eq:si1c} and \cref{eq:si1125c} to machine precision can be found at \url{https://github.com/AZaitzeff/SIgradflow}.
 
\cref{alg:m3rdorder} details our third order version for solving gradient flows with solution dependent inner product. The method adds another condition for stability to hold, namely:
\begin{equation}
\label{eq:opreq}
    \op(u)-\frac{1}{72}k^2D^2\op(u)(w,w)
\end{equation} needs to be positive definite for all $u$ and $w$. Now we will prove that \cref{alg:m3rdorder} produces a third order approximation.

By applying \cref{eq:rec}, the coefficients \cref{eq:si1c} give the following expansion for $u_{*_1}$:
\small
\begin{align}
\begin{split}
\label{eq:si1}
    u_{*_1}&=u_n-\frac{1}{6}k\op_n DE_n+\frac{1}{36}k^2\op_n D^2E_n(\op_n DE_n)+O(k^3)
    \end{split}
\end{align}
\normalsize
Now we can expand $\bar{u}$ by using \cref{eq:rec} and expanding $u_{*_1}$ around $u_n$
\small
\begin{align}
\begin{split}
\label{eq:3o1h}
    \bar{u}&=u_n-\frac{1}{2}k\op(u_{*_1}) DE_n+\frac{1}{8}k^2\op(u_{*_1})D^2E_n(\op(u_{*_1})DE_n)\\&-\frac{1}{48}k^3 \op(u_{*_1})D^2E_n\left(\op(u_{*_1})D^2E_n\left(\op(u_{*_1})DE_n\right)\right)\\&-\frac{1}{48}k^3\op(u_{*_1})D^3E_n\big(\op(u_{*_1})DE_n,\op(u_*)DE_n\big)\\&+\frac{1}{144}k^3D^2\op(u_{*_1})\big(\op(u_{*_1})DE(u_{*_1}),\op(u_{*_1})DE(u_{*_1})\big)DE_n+O(k^4)\\
    &=u_n-\frac{1}{2}k\op_n DE_n+\frac{1}{12}k^2D\op_n(\op_nDE_n) DE_n+\frac{1}{8}k^2\op_nD^2E_n(\op_nDE_n)\\&-\frac{1}{72}k^3D\op_n(\op_n D^2E_n(\op_n DE_n)) DE_n\\&-\frac{1}{48}k^3\op_nD^2E_n(D\op_n(\op_nDE_n)DE_n)\\&-\frac{1}{48}k^3D\op_n(\op_nDE_n) D^2E_n(\op_n DE_n)\\&-\frac{1}{48}k^3 \op_n D^2E_n\left(\op_n D^2E_n\left(\op_n DE_n\right)\right)\\&-\frac{1}{48}k^3\op_n D^3E_n\big(\op_n DE_n,\op_n DE_n\big)+O(k^4)
    \end{split}
\end{align}
\normalsize
Now we will apply the same steps to derive the expansions of $u_{*_{2,1}}$
\small
\begin{align}
\begin{split}
\label{eq:1s}
    u_{*_{2,1}}&=u_n-\frac{2}{5}k\op_n DE_n+O(k^2)
    \end{split}
\end{align}
\normalsize
and $u_{*_{2,2}}$
\small
\begin{align}
\begin{split}
\label{eq:si1125}
    u_{*_{2,2}}&=u_n-\frac{5}{6}k\op(u_{*_{2,1}}) DE_n+\frac{11}{36}k^2\op(u_{*_{2,1}}) D^2E_n(\op(u_{*_{2,1}}) DE_n)+O(k^3)\\
    &=u_n-\frac{5}{6}k\op_n DE_n\\&+\frac{1}{3}k^2D\op_n(\op_n DE_n) DE_n+\frac{11}{36}k^2\op_n D^2E_n(\op_n DE_n)+O(k^3)\\
    \end{split}
\end{align}
\normalsize
Finally, we can find the expansion of $u_{n+1}$. We will first apply \cref{eq:rec} around $\bar{u}$
\small
\begin{align*}
u_{n+1}&=\bar{u}-\frac{1}{2}k\op(u_{*_{2,2}}) DE(\bar{u})+\frac{1}{8}k^2\op(u_{*_{2,2}})D^2E(\bar{u})(\op(u_{*_{2,2}})DE(\bar{u}))\\&-\frac{1}{48}k^3 \op(u_{*_{2,2}})D^2E(\bar{u})\left(\op(u_{*_{2,2}})D^2E(\bar{u})\left(\op(u_{*_{2,2}})DE(\bar{u})\right)\right)\\&-\frac{1}{48}k^3\op(u_{*_{2,2}})D^3E(\bar{u})\big(\op(u_{*_{2,2}})DE(\bar{u}),\op(u_{*_{2,2}})DE(\bar{u})\big)\\&+\frac{1}{144}k^3D^2\op(u_{*_{2,2}})\big(\op(u_{*_{2,2}})DE(u_{*_{2,2}}),\op(u_{*_{2,2}})DE(u_{*_{2,2}})\big)DE(\bar{u})+O(k^4)
\end{align*}
expand $u_{*_{2,2}}$
\begin{align*}  
    u_{n+1}&=\bar{u}-\frac{1}{2}k\op_n DE(\bar{u})+\frac{5}{12}k^2D\op_n(\op_n DE_n) DE(\bar{u})+\frac{1}{8}k^2\op_n D^2E(\bar{u})(\op_n DE(\bar{u}))\\&-\frac{1}{6}k^3D\op_n(D\op_n(\op_n DE_n) DE_n) DE(\bar{u})\\&-\frac{1}{6}k^3D^2\op_n(\op_n DE_n,\op_n DE_n) DE_n) DE(\bar{u})\\&-\frac{11}{72}k^3D\op_n(\op_n D^2E_n(\op_n DE_n)) DE(\bar{u})\\&
    -\frac{5}{48}k^3 D\op_n(\op_n DE(\bar{u})) D^2E(\bar{u})(\op_n DE(\bar{u}))\\&-\frac{5}{48}k^3\op_n D^2E(\bar{u})(D\op_n(\op_n DE(\bar{u})) DE(\bar{u}))\\&-\frac{1}{48}k^3 \op_n D^2E(\bar{u})\left(\op_n D^2E(\bar{u})\left(\op_n DE(\bar{u})\right)\right)\\&-\frac{1}{48}k^3\op_nD^3E(\bar{u})\big(\op_n DE(\bar{u}),\op_n DE(\bar{u})\big)+O(k^4)
    \end{align*}
then expand $\bar{u}$ around $u_n$:
    \begin{align*}
    u_{n+1}&=u_n-k\op_n DE_n+\frac{1}{2}k^2D\op_n(\op_n DE_n) DE_n+\frac{1}{2}k^2\op_n D^2E_n(\op_n DE_n)\\&-\frac{1}{6}k^3D\op_n(D\op_n(\op_n DE_n) DE_n) DE_n\\&-\frac{1}{6}k^3D^2\op_n(\op_n DE_n,\op_n DE_n) DE_n) DE_n\\&-\frac{1}{6}k^3D\op_n(\op_n D^2E_n(\op_n DE_n)) DE_n\\&
    -\frac{1}{3}k^3 D\op_n(\op_n DE_n) D^2E_n(\op_n DE_n)\\&-\frac{1}{6}k^3\op_n D^2E_n(D\op_n(\op_n DE_n) DE_n)\\&-\frac{1}{6}k^3 \op_n D^2E_n\left(\op_n D^2E_n\left(\op_n DE_n\right)\right)\\&-\frac{1}{6}k^3\op_nD^3E_n\big(\op_n DE_n,\op_n DE_n\big)+O(k^4)\\
\end{align*}
\normalsize
The Taylor expansion of $u_{n+1}$ matches \cref{eq:truewmu} to third order. As long as \cref{eq:opreq} holds, \[E(u_{n+1})\leq E(\bar{u})\leq E(u_n)\]
by \cref{claim:ms}.

\begin{remark}
\label{remark:fullyimp}
In \cref{alg:m2ndorder} and \cref{alg:m3rdorder} we can instead handle $E(u)$ fully implicitly as we do in \cite{alexander2019variational}. We need to substitute higher order implicit methods for the corresponding semi-implicit methods. We give theses fully implicit versions in \cref{sec:append}.
\end{remark}

\section{Numerical Examples}
\label{sec:nr}
In this section, we will apply the second and third order accurate conditionally stable schemes to a variety of gradient flows, some with fixed inner product and some with solution dependent inner product.
Careful numerical convergence studies are presented in each case to verify the anticipated convergence rates of previous sections.

\subsection{Gradient Flows with Fixed Inner Product}

\begin{table}
\begin{center}
\begin{tabular}{|c|c|c|c|c|c|}
\hline
Number of& & & & &\\time steps&$2^5$&$2^{6}$&$2^{7}$&$2^{8}$&$2^{9}$\\
\hline
$L^2$ error (2nd order)&5.28e-05&1.16e-05& 2.71e-06&6.58e-07&1.62e-07\\
\hline
Order&-&2.19&2.09&2.04 &2.02 \\
\hline
$L^2$ error (3rd order)&1.11e-06&7.44e-07& 1.62e-07&2.51e-08&3.36e-09\\
\hline
Order&-&0.57&2.21&2.68 &2.90 \\
\hline
\end{tabular}
\caption{\footnotesize The new second and third order accurate, unconditionally stable schemes (see \cref{remark:fullyimp}) for gradient flows porous medium equation.}
\label{tab:PMEsi}
\end{center}
\end{table}

\begin{figure}[h]
\begin{center}
\includegraphics[width=.45\textwidth]{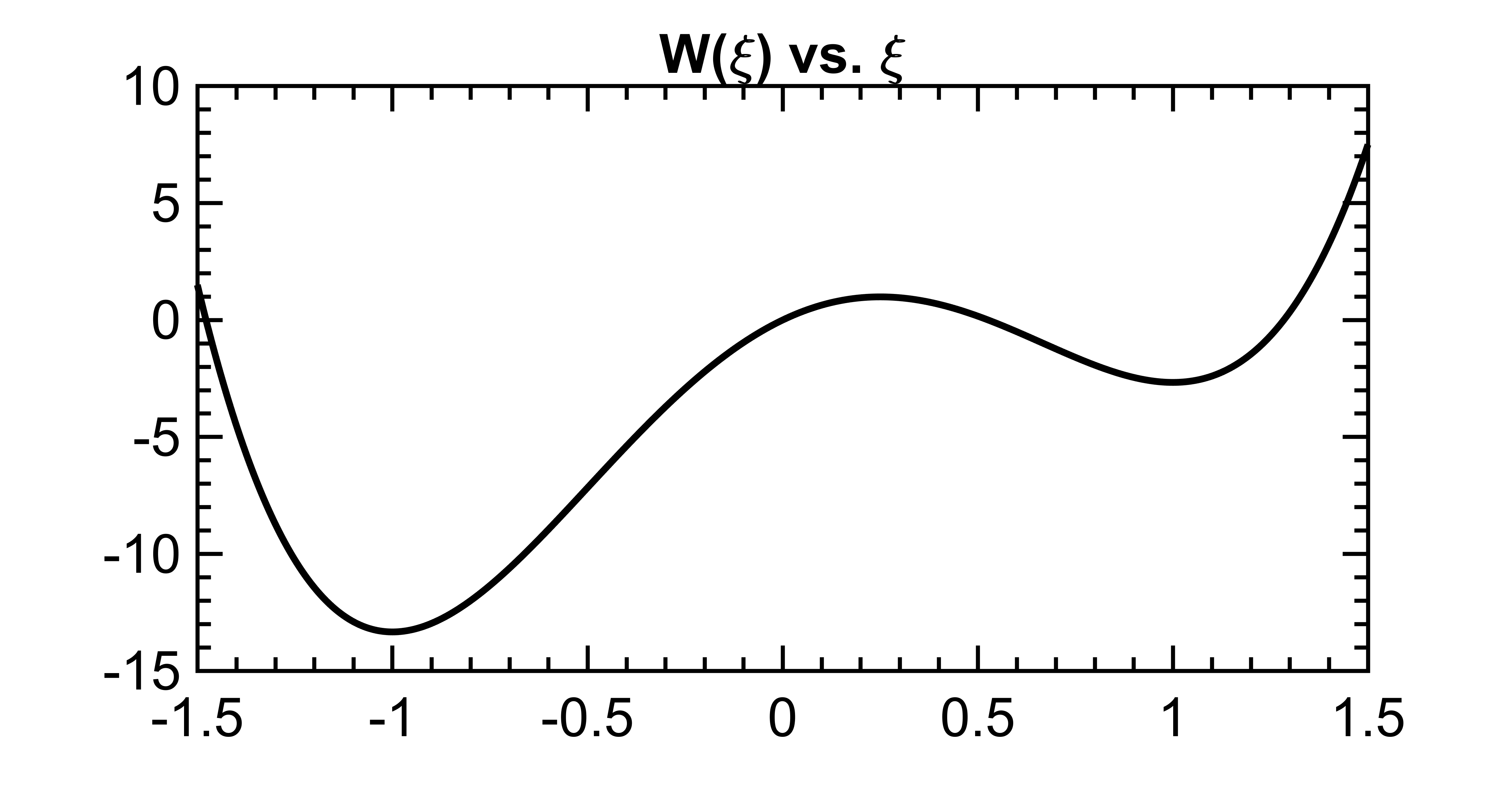}
\includegraphics[width=.45\textwidth]{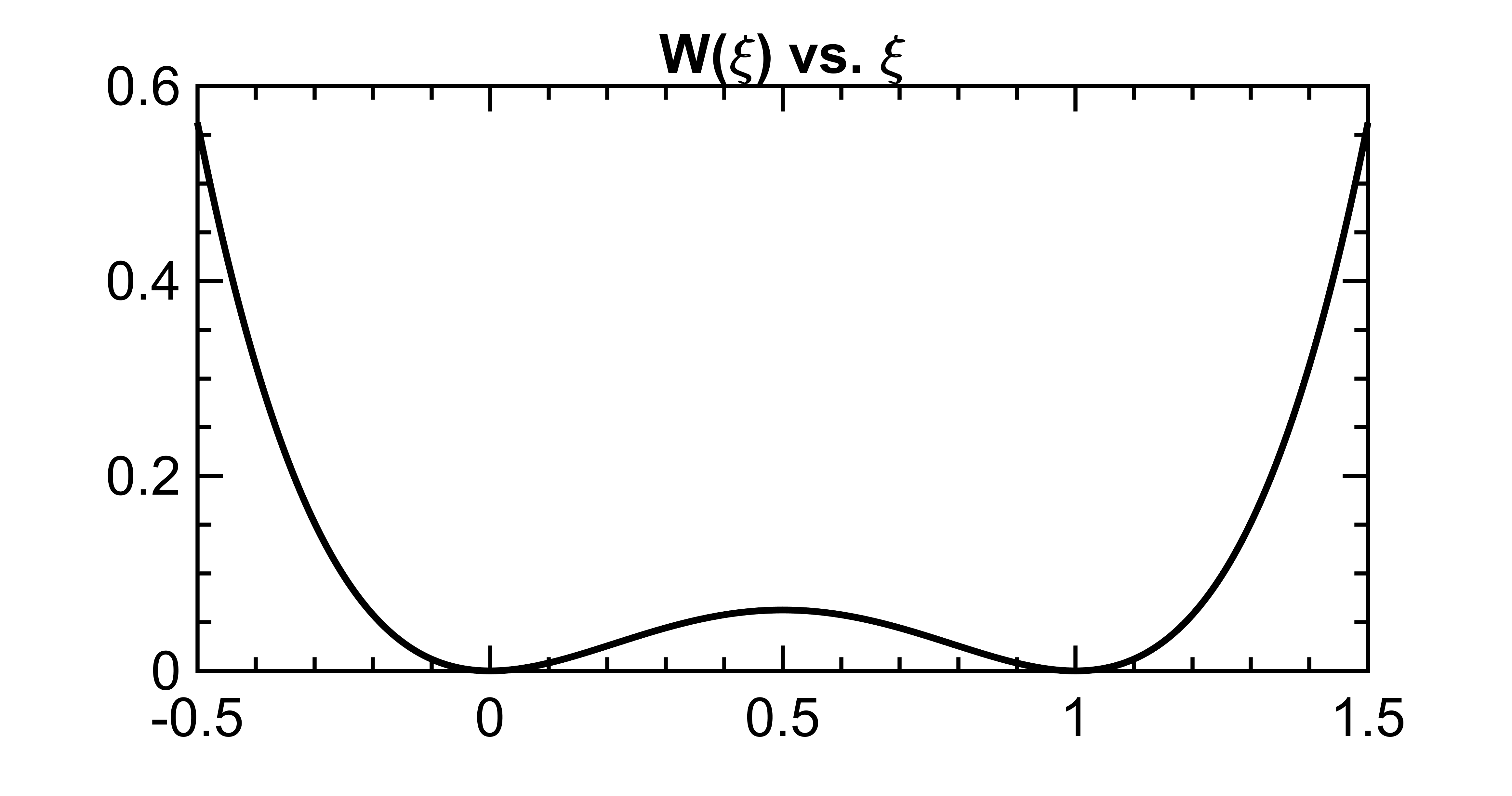}
\caption{\footnotesize The double well potentials used in the Allen-Cahn \cref{eq:ac} and Cahn-Hilliard \cref{eq:ch} equations: One with unequal depth wells and the other with equal depth wells.}
\label{fig:W}
\end{center}
\end{figure}

\begin{figure}[h]
  \begin{center}
\includegraphics[width=.5\textwidth]{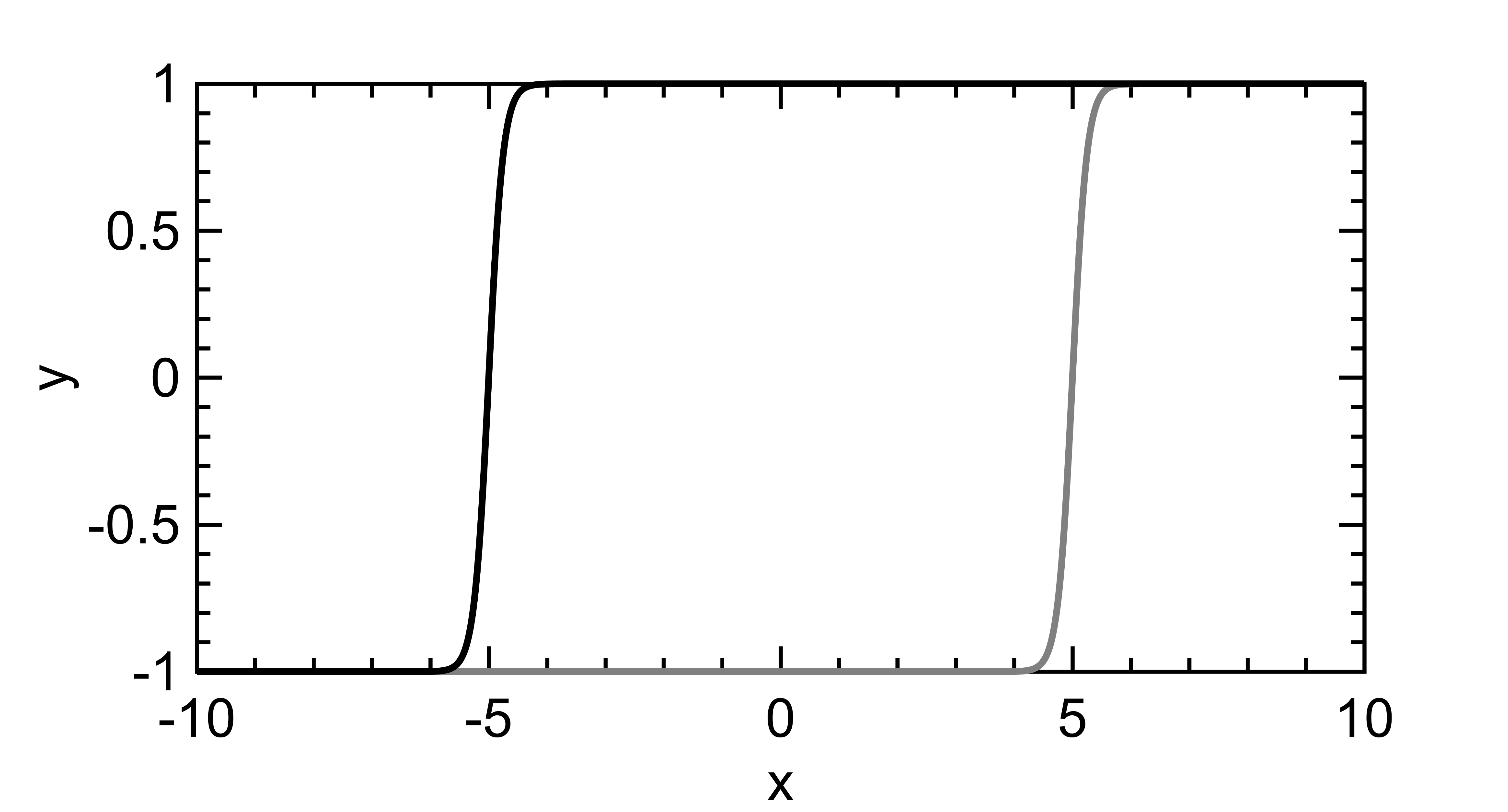} %
\end{center}
\caption{\footnotesize The initial condition (black) and the solution at final time (gray) in the numerical convergence study on the 1D Allen-Cahn equation \cref{eq:ac} with a potential that has unequal depth wells.}
\label{fig:1dac}
\end{figure}

We start with the Allen-Cahn equation
\begin{equation}
\label{eq:ac}
u_t = \Delta u - W'(u)
\end{equation}
where $W:\mathbb{R}\to\mathbb{R}$ is a double-well potential.
This corresponds to gradient flow for the energy
\begin{equation}
\label{eq:acE}
E(u) = \int \frac{1}{2} \|\nabla u\|^2 + W(u) \, dx
\end{equation}
with respect to the $L^2$ inner product.

First, we consider equation \cref{eq:ac} in one space dimension, with the potential $W(u) = 8u-16u^2-\frac{8}{3}u^3+8u^4$.
This is a double well potential with unequal depth wells; see \cref{fig:W}.
In this case, equation \cref{eq:ac} is well-known to possess traveling wave solutions on $x\in\mathbb{R}$, see \cref{fig:1dac}.
We choose the initial condition $u(x,0) = \tanh(4x + 20)$; the exact solution is then $u_*(x,t) = \tanh(4x + 20 - 8t)$.
The computational domain is $x\in[-10,10]$, discretized into a uniform grid of $8193$ points.
We approximate the solution on $\mathbb{R}$ by using the Dirichlet boundary conditions $u(\pm 10,t) = \pm 1$: The domain size is large enough that the mismatch in boundary conditions do not substantially contribute to the error in the approximate solution over the time interval $t\in[0,5]$. We use $E_1(u)=\int \frac{1}{2}|\nabla u|^2 dx$ and $E_2(u)=\int W(u) dx$.
\Cref{tab:twms5} tabulates the error in the computed solution at time $T=5$ for our two new schemes. 

\begin{table}
\begin{center}
\begin{tabular}{|c|c|c|c|c|c|}
\hline
Number of& & & & &\\time steps&$2^9$&$2^{10}$&$2^{11}$&$2^{12}$&$2^{13}$\\
\hline
$L^2$ error (2nd order) &2.08e-01&5.96e-02&1.61e-02& 4.22e-03&1.08e-03\\
\hline
Order&-&1.81&1.89&1.94 &1.97  \\
\hline
$L^2$ error (3rd order) &2.06e-03&3.26e-04& 4.68e-05&6.32e-06&8.33e-07\\
\hline
Order&-&2.66&2.80&2.89 &2.92  \\
\hline
\end{tabular}
\caption{\footnotesize The new second \cref{eq:2ndordergamma} and third \cref{eq:3rdordergamma} order accurate, conditionally stable schemes \cref{eq:ms} on the one-dimensional Allen-Cahn equation \cref{eq:ac} with a traveling wave solution.}
\label{tab:twms5}
\end{center}
\end{table}

\begin{figure}[h]
  \begin{center}
\includegraphics[width=.45\textwidth]{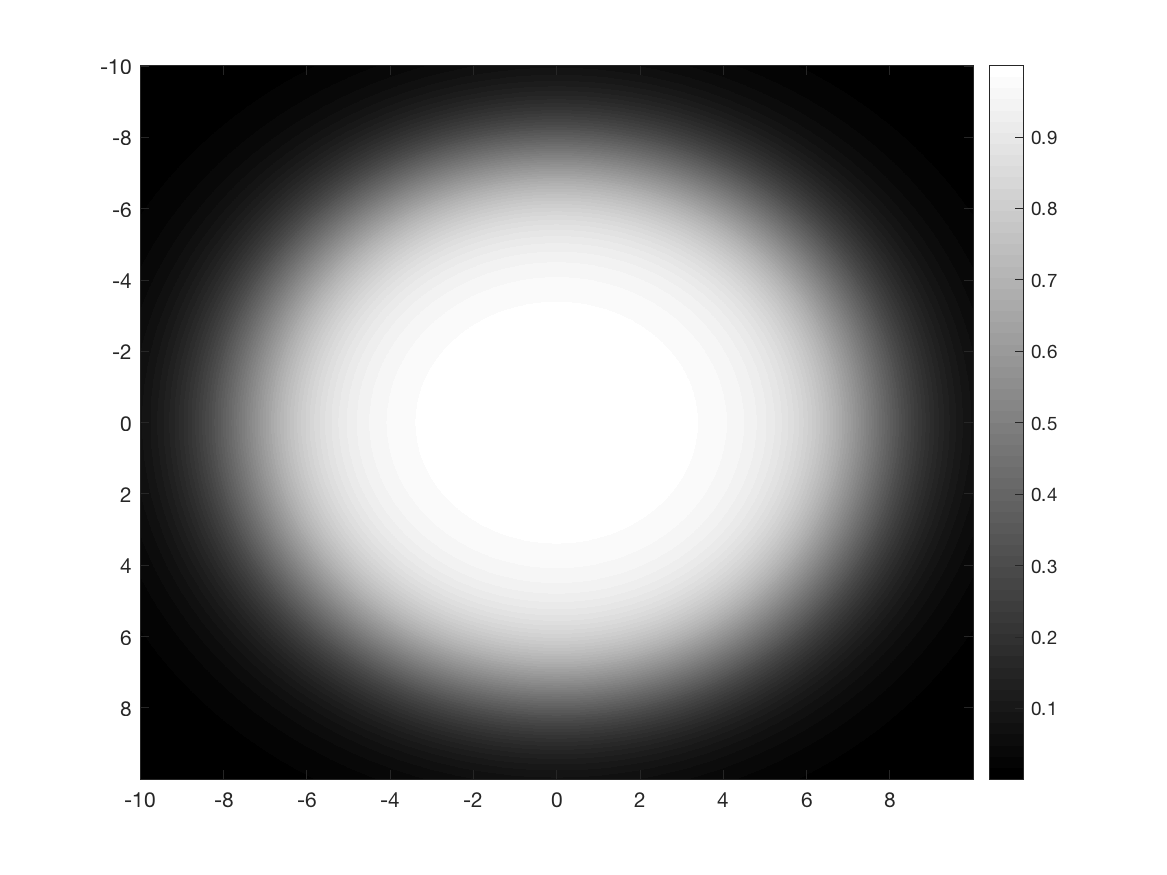}
\includegraphics[width=.45\textwidth]{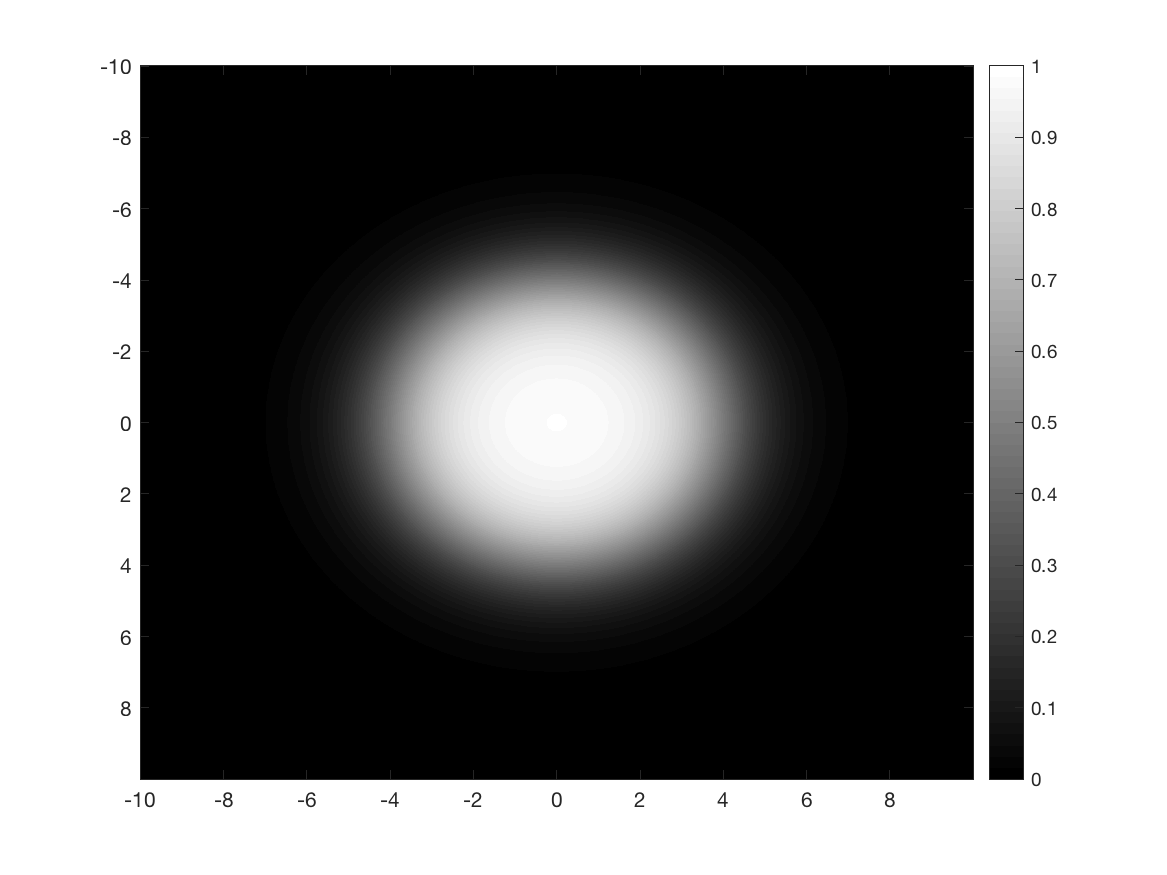}
\end{center}
\caption{\footnotesize Initial condition and the solution at final time for the 2D Allen-Cahn equation with a potential that has equal depth wells.}
\label{fig:2dac}
\end{figure}

Next, we consider the Allen-Cahn equation \cref{eq:ac} in two space-dimensions, with the potential $W(u) = u^2(1-u)^2$ that has equal depth wells; see \cref{fig:W}.
We take the initial condition $u(x,y,0)=\frac{1}{1+\exp[-(7.5-\sqrt{x^2+y^2})]}$ on the domain $x\in [-10,10]^2$, and impose periodic boundary conditions. Once again we use $E_1(u)=\int \frac{1}{2}\|\nabla u\|^2 dx$ and $E_2(u)=\int W(u) dx$.
As a proxy for the exact solution of the equation with this initial data, we compute a very highly accurate numerical approximation $u_*(x,y,t)$ via the following second order accurate in time, semi-implicit, multi-step scheme ~\cite{chen1998applications} on an extremely fine spatial grid and take very small time steps:
\begin{equation*}
    \frac{3}{2}u^{n+1}-2u^{n}+\frac{1}{2}u^{n-1}=k\Delta u^{n+1}-k(2W'(u^{n})-W'(u^{n-1})).
\end{equation*}
\Cref{tab:acms62d} show the errors and convergence rates for the approximate solutions computed by our new multi-stage schemes.

\begin{table}[h]
\begin{center}
\begin{tabular}{|c|c|c|c|c|c|}
\hline
Number of& & & & &\\time steps&$2^8$&$2^9$&$2^{10}$&$2^{11}$&$2^{12}$\\
\hline
$L^2$ error (2nd order)& 3.62e-05&9.07e-06&2.27e-06& 5.68e-07 &1.41e-07\\
\hline
Order&-&2.00&2.00&2.00&2.00 \\
\hline
$L^2$ error (3rd order)&2.35e-05&3.18e-06& 4.15e-07&5.29e-08&6.24e-09\\
\hline
Order&-&2.88&2.94&2.97 &3.08 \\
\hline
\end{tabular}
\caption{\footnotesize The new second \cref{eq:2ndordergamma} and third \cref{eq:3rdordergamma} order accurate, conditionally stable schemes \cref{eq:ms} on the two-dimensional Allen-Cahn equation \cref{eq:ac} with a potential that has equal depth wells.}
\label{tab:acms62d}
\end{center}
\end{table}

\begin{figure}[h]
  \begin{center}
\includegraphics[width=.45\textwidth]{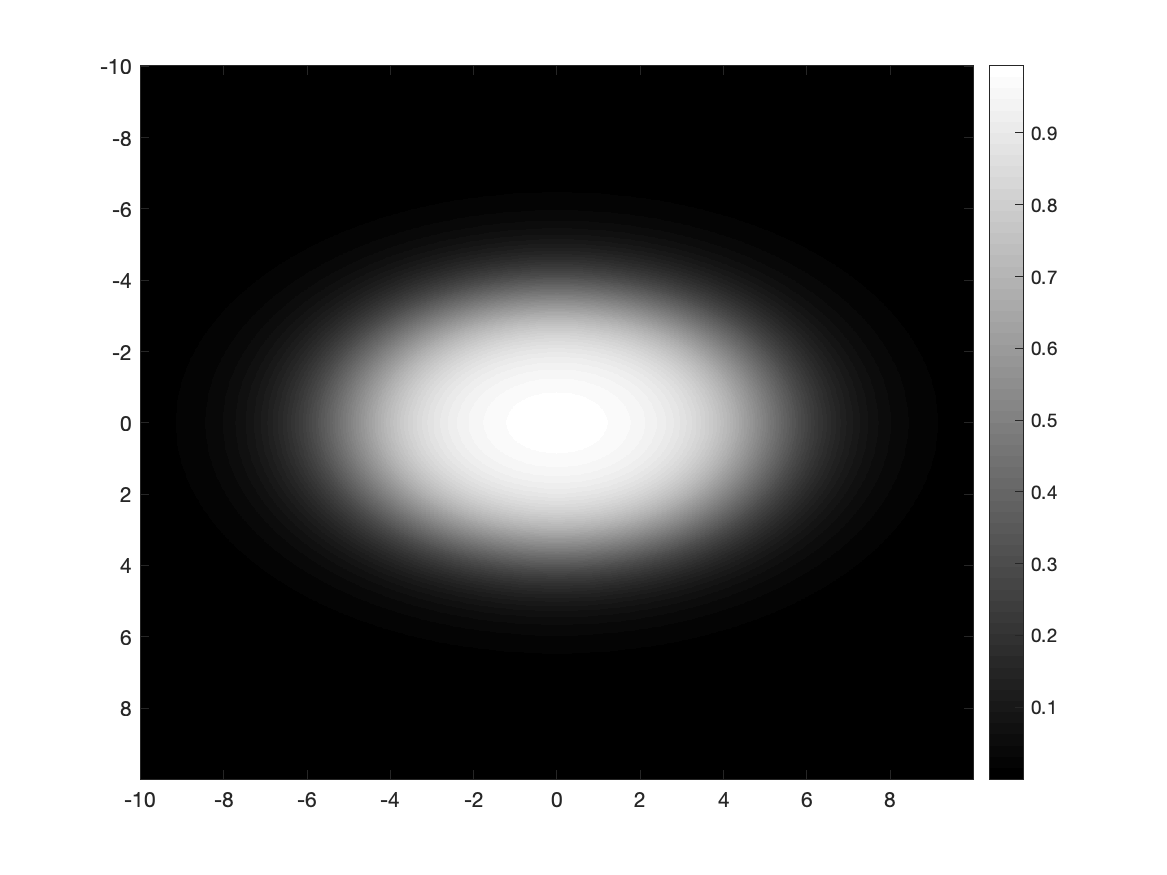}
\includegraphics[width=.45\textwidth]{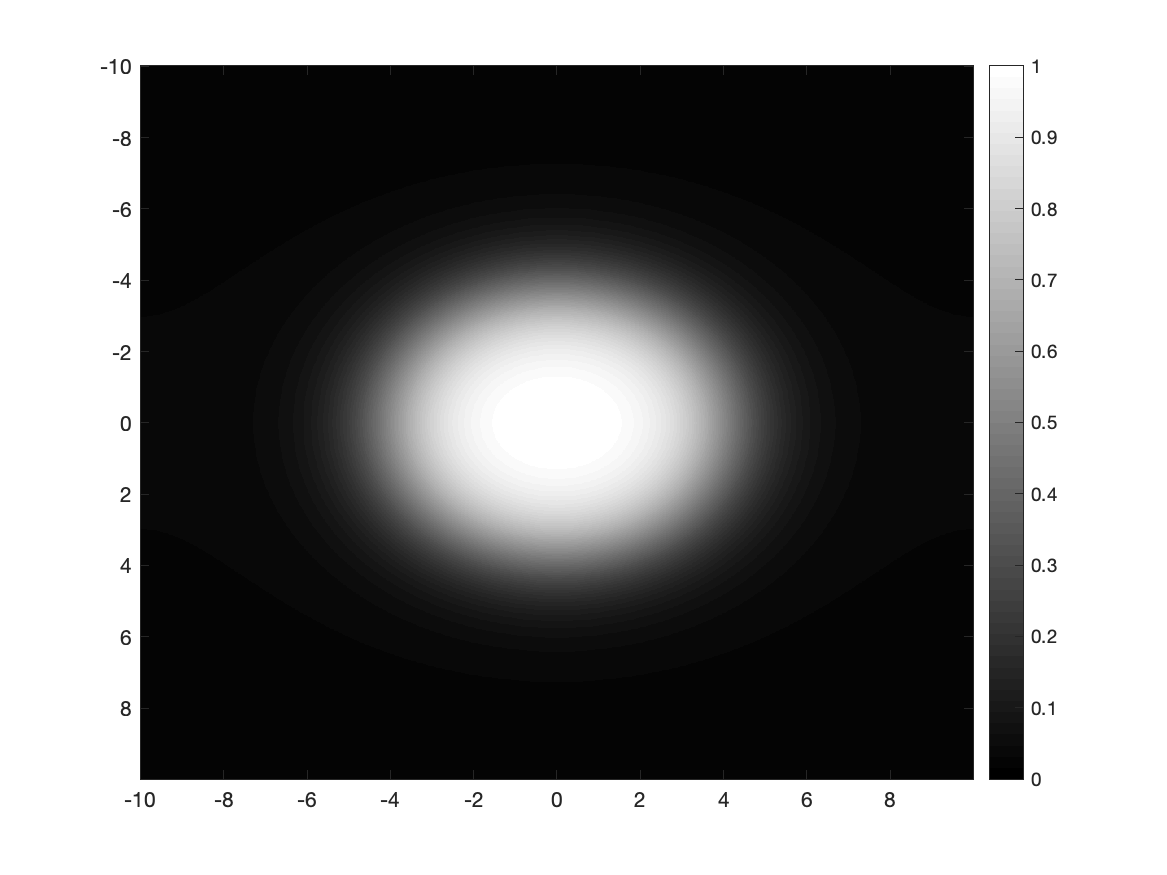}
\end{center}
\caption{\footnotesize Initial condition and the solution at final time for the 2D Cahn-Hillard equation with a potential that has equal depth wells.}
\label{fig:2dch}
\end{figure}

For our next example, we consider the Cahn-Hilliard equation
\begin{equation}
\label{eq:ch}
u_t = -\Delta \big( \Delta u - W'(u) \big)
\end{equation}
where we take $W$ to be the double well potential $W(u) = u^2(1-u)^2$ with equal depth wells and impose periodic boundary conditions.
This flow is also gradient descent for energy \cref{eq:acE}, but with respect to the $H^{-1}$ inner product:
\begin{equation*}
\langle u \,, v \, \rangle = \int u \Delta^{-1} v \, dx.
\end{equation*}
Starting from the initial condition $u(x,y,0)=\frac{1}{1+\exp[-(5-\sqrt{x^2+y^2})]}$, we computed a proxy for the ``exact'' solution once again using the second order accurate, semi-implicit multi-step scheme from \cite{chen1998applications}:
\begin{equation*}
    \frac{3}{2}u^{n+1}-2u^{n}+\frac{1}{2}u^{n-1}=-k\Delta[\Delta u^{n+1}-(2W'(u^{n})-W'(u^{n-1}))]
\end{equation*}
where the spatial and temporal resolution was taken to be high to ensure the errors are small.
\Cref{tab:chms32d} show the errors and convergence rates for the approximate solutions computed by our new multi-stage schemes.

\begin{table}
\begin{center}
\begin{tabular}{|c|c|c|c|c|c|}
\hline
Number of& & & & &\\time steps&$2^{7}$&$2^{8}$&$2^{9}$&$2^{10}$&$2^{11}$\\
\hline
$L^2$ error (2nd order)& 6.20e-04 & 1.92e-04  & 5.59e-05& 1.55e-05&4.09e-06\\
\hline
Order&-&1.69&1.78&1.85&1.92 \\
\hline
$L^2$ error (3rd order) &6.45e-06&1.35e-06& 2.51e-07&4.15e-08&7.20e-09\\
\hline
Order&-&2.25&2.43&2.60&2.53 \\
\hline
\end{tabular}
\caption{\footnotesize The new second \cref{eq:2ndordergamma} and third \cref{eq:3rdordergamma} order accurate, conditionally stable schemes \cref{eq:ms} on the two-dimensional Cahn-Hilliard equation \cref{eq:ch} with a potential that has equal depth wells.}
\label{tab:chms32d}
\end{center}
\end{table}

\begin{figure}[h]
  \begin{center}
\includegraphics[width=.5\textwidth]{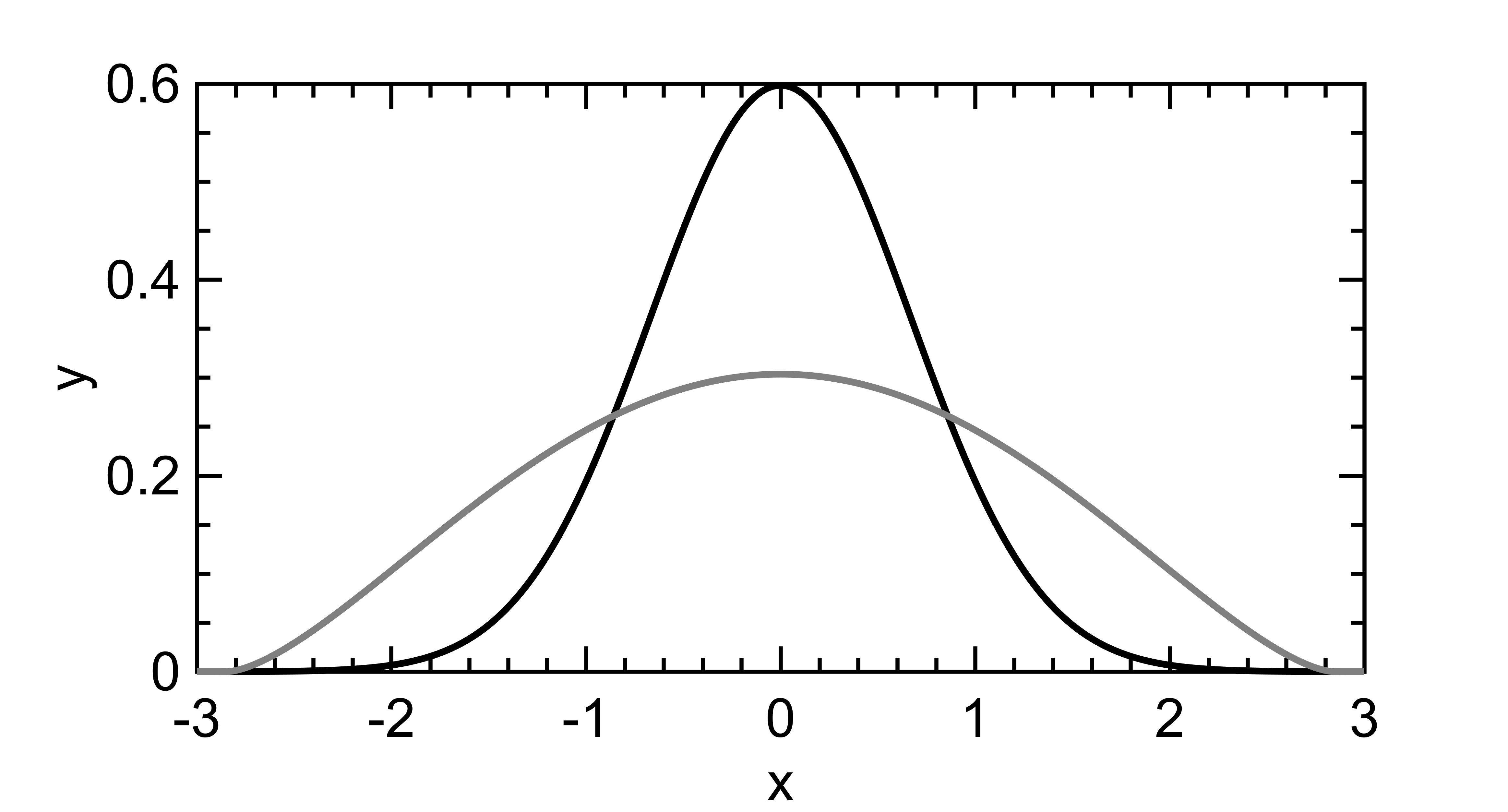} %
\end{center}
\caption{\footnotesize The initial condition (black) and the solution at final time (gray) in the numerical convergence study on the PME}
\label{fig:1dpme}
\end{figure}

As a final example we do the following porous medium equation:
\begin{equation}
\label{eq:pme53}
    u_t=\Delta u^{5/3}
\end{equation}

Under the $H^{-1}$ inner product, \cref{eq:pme53} is gradient flow for the energy 
\[
E(u)=\frac{3}{8} \int u^{8/3} dx.
\]
Our initial data is
\begin{equation}
\label{eq:pmeinit}
    u(x,0)=\frac{3}{2\sqrt{2\pi}}\exp\bigg(-\frac{9x^2}{8}\bigg)
\end{equation}
in $x\in [-3,3]$ with derivative zero Neumann boundary conditions. We run the simulation for $T=1$. See \cref{fig:1dpme} for our initial and final curve. We generate the ``true'' solution using the L-stable (but not energy stable) TR-BDF2 method with a high spatial and temporal resolution. See \cref{tab:pme1} for results.

\begin{table}
\begin{center}
\begin{tabular}{|c|c|c|c|c|c|}
\hline
Number of& & & & &\\time steps&$2^{12}$&$2^{13}$&$2^{14}$&$2^{15}$&$2^{16}$\\
\hline
$L^2$ error (2nd order)& 1.91e-06 & 6.85e-07  & 2.26e-07& 6.90e-08&1.97e-08\\
\hline
Order&-&1.48&1.60&1.71&1.81 \\
\hline
$L^2$ error (3rd order) &2.21e-07&5.69e-08& 1.26e-08&2.41e-09&4.13e-10\\
\hline
Order&-&1.96&2.18&2.38&2.54 \\
\hline
\end{tabular}
\caption{\footnotesize The new second \cref{eq:2ndordergamma} and third \cref{eq:3rdordergamma} order accurate, conditionally stable schemes \cref{eq:ms} on the porous medium equation}
\label{tab:pme1}
\end{center}
\end{table}

\subsection{Gradient Flow For Solution Dependent Inner Product}
Our first example we present in this section is the heat equation, $u_t=\Delta u$, but with a different energy. Under the Wasserstein metric (denoted as $W_2$), the heat equation is a gradient flow for the negative entropy \cite{jordan1998variational}:
\begin{equation}
   \label{SIeq:entropy}
   E(u)=\int u\log(u) dx. 
\end{equation}
However the minimization \[
\argmin_u E(u)+\frac{1}{2k}W_2^2(u,u_n)
\]
is a difficult optimization problem. On the other hand, we can approximate the the Wasserstein metric, $W_2(u,v)$, with
\begin{equation}
\label{eq:WIP}
\langle u-v,\op(u)^{-1}(u-v) \rangle_{L^2}\text{ where }\op(u)=-\nabla \cdot u \nabla
\end{equation}

when $u$ and $v$ are near each other.
Indeed 
\[-\op(u)\nabla_{L^2} E(u)=\nabla \cdot u \nabla(\log(u)+1)=\Delta u.\]
Thus, we can alternatively think of the heat equation as minimizing movements on negative entropy with respect to the solution dependent inner product \cref{eq:WIP} and therefore use \cref{alg:m2ndorder} and \cref{alg:m3rdorder} to evolve the heat equation while decreasing the negative entropy \cref{SIeq:entropy} at every step. 

We use the exact solution 
$u(x,t)=\cos(\pi x)\exp(-t\pi^2)+2$ as our test with domain $x\in [0,1]$ using derivative zero Neumann boundary conditions. Our initial data is $u(x,0)$ and we run the simulation to final time $T=\frac{1}{10}$. We use $E_1(u)=\frac{1}{2} \int u^2 dx$ and $E_2(u)=\int u\log(u) dx-\frac{1}{2}\int u^2 dx$ in \cref{eq:ms} so at every step we are solving a linear systems of equation. 
 We run simulation for $T=\frac{1}{10}$. See \cref{tab:heat} for results.

\begin{table}
\begin{center}
\begin{tabular}{|c|c|c|c|c|c|}
\hline
Number of& & & & &\\time steps&$2^3$&$2^{4}$&$2^{5}$&$2^{6}$&$2^{7}$\\
\hline
$L^2$ error (2nd order)&1.06e-03&3.11e-04&8.58e-05& 2.27e-05&5.85e-06\\
\hline
Order&-&1.77&1.86&1.92&1.96   \\
\hline
$L^2$ error (3rd order)&1.00e-05&1.57e-06& 2.20e-07&3.04e-08&4.29e-09\\
\hline
Order&-&2.69&2.82&2.87 &2.83  \\
\hline
\end{tabular}
\caption{\footnotesize The new second (\cref{alg:m2ndorder}) and third (\cref{alg:m3rdorder}) order accurate, conditionally stable schemes for gradient flows with solution dependent inner product on the heat equation with Wasserstein metric.}
\label{tab:heat}
\end{center}
\end{table}
The next example is the porous medium equation in one dimension. The energy is
\[
E(u)=\frac{3}{2}\int u^{5/3} dx
\]
under the Wasserstein metric. As with the heat equation, we can again replace the Wasserstein metric with \cref{eq:WIP}. We will let $E_1(u)=E(u)$ and $E_2(u)=0$. We use the same test as in the $H^{-1}$ gradient flow porous medium equation (see \cref{eq:pmeinit} and accompanying explanation). We present the results of the porous medium equation test with movement limiter \cref{eq:WIP} in \cref{tab:pme2}.

\begin{table}
\begin{center}
\begin{tabular}{|c|c|c|c|c|c|}
\hline
Number of& & & & &\\time steps&$2^{4}$&$2^5$&$2^{6}$&$2^{7}$&$2^{8}$\\
\hline
$L^2$ error (2nd order)& 2.69e-04& 6.10e-05& 1.49e-05&3.71e-06&9.25e-07 \\
\hline
Order&-&2.14&2.04&2.01 &2.00 \\
\hline
$L^2$ error (3rd order)&2.88e-0& 3.17e-06  & 3.72e-07  & 4.53e-08& 5.50e-09\\
\hline
Order&-&3.18&3.09&3.04 &3.04 \\
\hline
\end{tabular}
\caption{\footnotesize The new second and third order accurate, unconditionally stable schemes (see \cref{remark:fullyimp}) for gradient flows with solution dependent inner product on the porous medium equation with the linearized Wasserstein metric.}
\label{tab:pme2}
\end{center}
\end{table}

\begin{figure}[h]
  \begin{center}
\includegraphics[width=.6\textwidth]{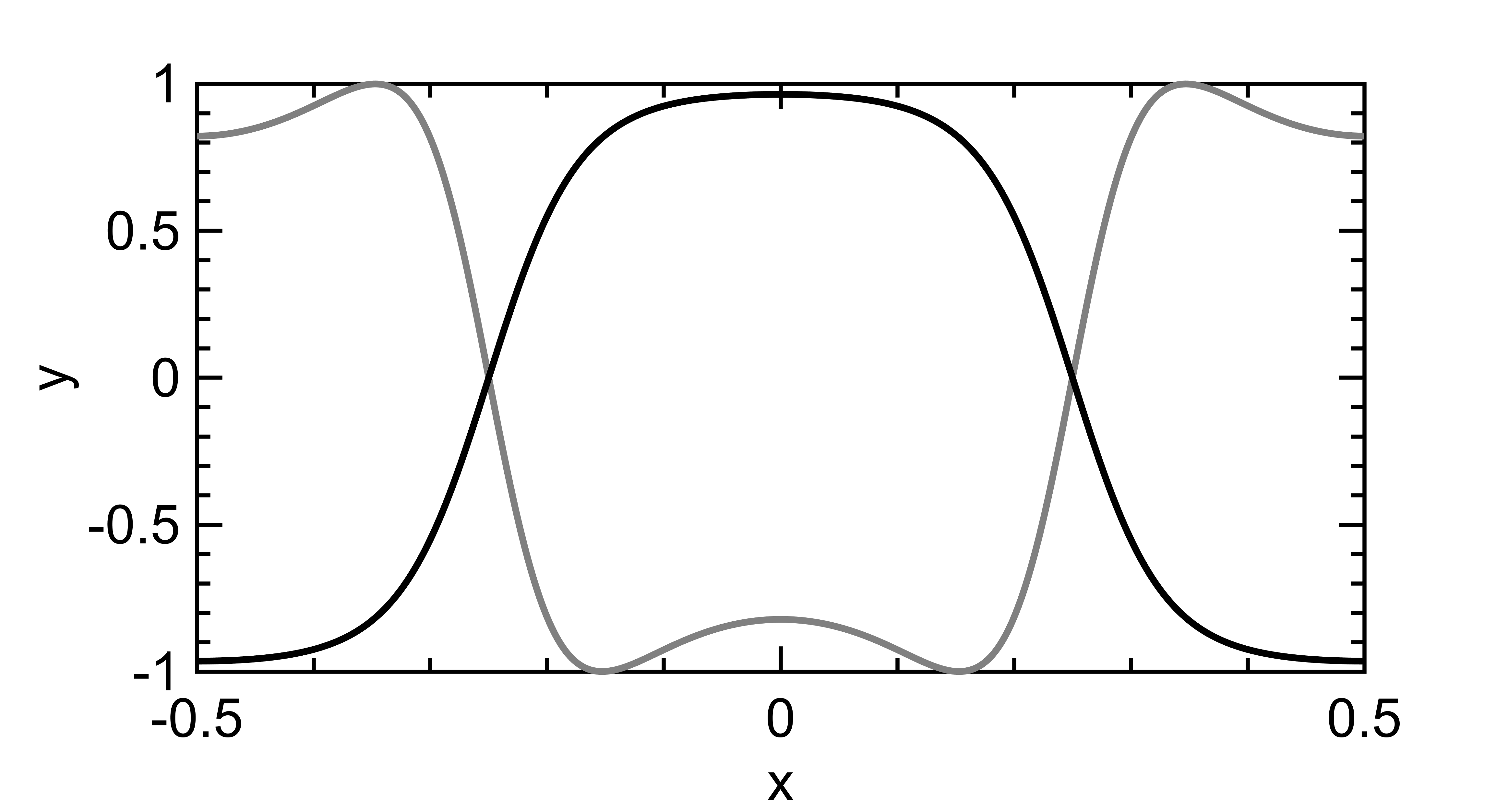}
\end{center}
\caption{\footnotesize Initial condition (black) and the solution at final time (gray) for the 1D Cahn-Hillard with variable mobility and forcing term example.}
\label{fig:1dchwm}
\end{figure}

For our final example, we consider the Cahn-Hilliard equation with variable mobility and a forcing term:
\small
\begin{equation}
\label{eq:chwm}
u_t = -\nabla \cdot \mu(u) \nabla \big( \epsilon^2 \Delta u - W'(u) -F(x)\big)
\end{equation}
\normalsize
where we take $W$ to be the double well potential $W(u) = (1-u^2)^2$ with equal depth wells, the forcing term to be $F(x)=\tanh\big(\frac{\cos(2\pi x)}{10\epsilon}\big)$ and the mobility to be $\mu(u) = (1-\epsilon)(1-u^2)^2+\epsilon$ to avoid degeneracy in the PDE.
This flow is gradient descent for energy 
\small
\begin{equation*}
E(u) = \int \frac{\epsilon^2}{2} \|\nabla u\|^2 + W(u)+uF(x) \, dx,
\end{equation*}
\normalsize
with respect to the solution dependent inner product
\small
\begin{equation}
\label{eq:CHWMIP}
\langle u-v,\op(u)^{-1}(u-v) \rangle_{L^2}\text{ where }\op(u)=-\nabla \cdot \mu{(u)} \nabla.
\end{equation}
\normalsize
For our example, we take $\epsilon=\frac{1}{20}$ and starting from the initial condition
\small \[u(x,0)=\tanh\bigg(\frac{\cos(2\pi x)}{10\epsilon}\bigg)\]
\normalsize
on the domain $x \in [-\frac{1}{2},\frac{1}{2}]^2$, and impose periodic boundary conditions. We run the PDE until time $T=\frac{1}{8}$.  We computed a proxy for the ``exact'' solution using the following second order BDF/AB scheme:
\small
\begin{multline*}
    3u^{n+1}+2k\epsilon^2\Delta^2 u^{n+1}=4u^n-u^{n-1}+4\big(k\epsilon^2\Delta^2 u^{n}-k\nabla\cdot \mu(u^n)\nabla[\epsilon^2\Delta u^{n}-W'(u^{n})]\big)\\-2\big(k\epsilon^2\Delta^2 u^{n-1}-k\nabla\cdot \mu(u^{n-1})\nabla[\epsilon^2\Delta u^{n-1}-W'(u^{n-1})]\big)
\end{multline*}
\normalsize
where the spatial and temporal resolution were taken to be high to ensure the errors are negligible. See \cref{fig:1dchwm} for plots of the initial condition and the solution at the final time.
\Cref{tab:chmswm} shows the errors and convergence rates for the approximate solutions computed by our new multi-stage schemes.

\begin{table}
\begin{center}
\begin{tabular}{|c|c|c|c|c|c|}
\hline
Number of& & & & &\\time steps&$2^{7}$&$2^{8}$&$2^{9}$&$2^{10}$&$2^{11}$\\
\hline
$L^2$ error (2nd order)& 4.10e-04 &  1.29e-04& 3.84e-05 &1.09e-05 & 2.96e-06 \\
\hline
Order&-&1.67&1.75&1.82&1.88 \\
\hline
$L^2$ error (3rd order) &1.79e-05& 3.68e-06&  6.72e-07 &  1.06e-07  & 1.42e-08\\
\hline
Order&-&2.28 &  2.46&  2.66 &  2.90 \\
\hline
\end{tabular}
\caption{\footnotesize  The new second \cref{alg:m2ndorder} and third \cref{alg:m3rdorder} order accurate, conditionally stable schemes for gradient flows with solution dependent inner product on the one-dimensional Cahn-Hilliard equation with variable mobility and forcing term \cref{eq:chwm}.}
\label{tab:chmswm}
\end{center}
\end{table}

 \section{Conclusion}
 We presented a new class of implicit-explicit additive Runge-Kutta schemes for gradient flows that are high order and conditionally stable. Additionally, we developed new high order stable schemes for gradient flows on solution dependent inner products. Both of these methods allow us to painlessly increase the order of accuracy of existing schemes for gradient flows without sacrificing stability. We provided many numerical examples of gradient flows, including those that have solution dependent inner product, and have shown that the methods achieve their advertised accuracy.
 
 However, in this paper, we have not developed a systematic approach to coming up with conditionally stable methods of a certain order. In fact, there may exist 2nd and 3rd order methods of fewer stages than given here. Additionally, whether these schemes can be used to achieve arbitrarily high (i.e. $\geq 4$) order in time is unknown. We leave these questions to future work.

 \appendix
 \section{A Second and Third Order Fully Implicit Methods for Gradient Flows with Solution Dependant Inner Product}
  \label{sec:append}
 Here we layout fully implicit, second and third order algorithms for gradient flows with solution dependant inner product,
\[
u_t=-\op(u)\nabla E(u).
\] 
In this case, each substep has form:
\begin{equation}
    \label{eq:Lfi}
    \bigg[\sum^{m-1}_{i=0}\gamma_{m,i}\bigg]U_{m}+k\op(u_*)\nabla_HE(U_m)=\sum^{m-1}_{i=0}\gamma_{m,i}U_i.
\end{equation}
The satisfaction of consistency to the exact solution, \cref{eq:truewmu}, is similar to its semi-implicit counterpart in \cref{sec:SDIP}. Theorem 2.1 in \cite{alexander2019variational} (which is the special case of this paper's \cref{claim:ms} when $E_2=0$) ensures that the coefficients in this section form multistage methods that are energy stable. The exact values for the coefficients for the multistage methods can be found at 
\url{https://github.com/AZaitzeff/SIgradflow}.
\subsection{Second order example} The following second order method is unconditionally energy stable.  
Fix a time step size $k>0$.  Set $u_n=u_0$. To obtain $u_{n+1}$ from $u_n$, carry out the following steps:
\begin{enumerate}

\item Find $u_*$ by solving
\small
\begin{equation*}
    u_*+\frac{1}{2}k\op_n\nabla E(u_*)=u_n
\end{equation*}
\normalsize
\item Find $u_{n+1}$ using \cref{eq:Lfi} with coefficients 
{\small
\begin{equation}
\label{eq:FI2ndordergamma}
    \gamma
\approx \left(
\begin{array}{ccc}
 5.0 & 0 & 0  \\
 -2.0 & 6.0 & 0  \\
 -2.0 & 0.22 & 6.29  \\
\end{array}\right).
\end{equation}}
and $u_*$ in $\op(u_*)$ given by Step 1 of this algorithm.
\end{enumerate}
\subsection{Third order example}
This third order algorithm is energy stable as long as
\begin{equation*}
    \op(u)-\frac{1}{72}k^2D^2\op(u)(w,w)
\end{equation*} is positive definite for all $u$ and $w$.

Fix a time step size $k>0$.
Set $u_n = u_0$. For convenience, we will denote
$D^2\op(u_*)\big(\op(u_*)\nabla E(u_*),\op(u_*)\nabla E(u_*)\big)$ as $D^2\op(u_*)$ and $MS(\tilde{k},\op(u_*),\tilde{u},\gamma)$ as $U_M$ in the multistage algorithm 
\small
\[
    \bigg[\sum^{m-1}_{i=0}\gamma_{m,i}\bigg]U_{m}+\tilde{k}\op(u_*)\nabla_HE(U_m)=\sum^{m-1}_{i=0}\gamma_{m,i}U_i.
\]
\normalsize
with time step $\tilde{k}$, operator $\op(u_*)$, $U_0=\tilde{u}$ and coefficients $\gamma$. To obtain $u_{n+1}$ from $u_n$, carry out the following steps:
\begin{enumerate}
\item Set $\gamma=\left(1\right)$ and $u_{*_1}=MS\Big(\frac{1}{6}k,\op(u_n),u_n,\gamma\Big)$.
\item Set 
\small
\begin{align}
\label{eq:FI3rdordergamma}
&\gamma \approx \left(
\begin{array}{cccccc}
 11.17& 0& 0&0&0&0  \\
 -7.5 & 19.43 & 0&0&0&0  \\
 -1.05& -4.75&13.98&0&0&0 \\
 1.8&0.05&-7.83&13.8&0& 0\\
 6.2&-7.17& -1.33&1.63&11.52&0\\
 -2.83&4.69&2.46&-11.55&6.68&11.95\\
\end{array}\right)
\end{align}
\normalsize
and \small\[\bar{u}=MS\bigg(\frac{1}{2}k,\op(u_{*_1})-\frac{1}{72}D^2\op(u_{*_1}),u_n,\gamma\bigg).\]\normalsize
\item Set $\gamma=\left(1\right)$ and $u_{*_{2,1}}=MS\Big(\frac{2}{5}k,\op(u_n),u_n,\gamma\Big)$.

\item Set
\small
\begin{equation}
\label{eq:FI1125}
\gamma \approx \left(
\begin{array}{cccc}
 6.17 & 0 & 0 & 0 \\
 -0.5 & 6 & 0 & 0 \\
 -3 & 2 & 7 & 0 \\
 -3.1 & 0 & 2.23 & 7.40 \\
\end{array}
\right)
\end{equation}
and \small\[
u_{*_{2,2}}=MS\bigg(\frac{5}{6}k,\op(u_{*_{2,1}}),u_n,\gamma\bigg).
\]\normalsize
\item Set $\gamma$ to \cref{eq:FI3rdordergamma} and then \small\[u_{n+1}=MS\bigg(\frac{1}{2}k,\op(u_{*_{2,2}})-\frac{1}{72}D^2\op(u_{*_{2,2}}),\bar{u},\gamma\bigg).\]\normalsize.
\end{enumerate}

\bibliographystyle{siamplain}
\bibliography{references}
\end{document}